\documentclass[11pt]{article}
\usepackage[margin = 2cm]{geometry}
\usepackage{amsmath,bbm,pdfsync}
\usepackage{amssymb}
\usepackage{amsthm}
\usepackage{xcolor}
\usepackage{csquotes}
\usepackage{enumerate}
\usepackage{subcaption}
\usepackage{todonotes}
 
\usepackage{graphicx,color}

\usepackage[russian, english]{babel}

\newtheorem{thm}{Theorem}[section]
\newtheorem{lem}[thm]{Lemma}
\newtheorem{cor}[thm]{Corollary}
\newtheorem{prop}[thm]{Proposition}

\DeclareMathOperator{\N}{\mathbb{N}}
\DeclareMathOperator{\R}{\mathbb{R}}

\newcommand{\iid}{i.i.d.\@ }
\allowdisplaybreaks

\usepackage[colorlinks=true,pdfpagemode=UseNone,urlcolor=blue,linkcolor=blue,citecolor=blue,breaklinks=true]{hyperref}

\numberwithin{equation}{section}

\newcommand{\vphi}{\varphi}
\newcommand{\cT}{{\cal T}}
\newcommand{\be}{\begin{equation}}
\newcommand{\ee}{\end{equation}}
\newcommand{\bal}{\begin{aligned}}
\newcommand{\enbal}{\end{aligned}}
\newcommand{\Pm}{{\mathbb P}}
\newcommand{\Rm}{{\mathbb R}}
\newcommand{\Nm}{{\mathbb N}}
\newcommand{\Em}{{\mathbb E}}
\newcommand{\Zm}{{\mathbb Z}}

\newcommand{\one}{{\mathbbm 1}}
\newcommand{\red}{\textcolor{red}}
\newcommand{\eps}{\varepsilon}
\newcommand{\farc}{\frac}
\newcommand{\vtheta}{\vartheta}
\newcommand{\med}{\hbox{med}}
 
\DeclareMathOperator{\supp}{\mathrm{supp}}

\newcommand{\commentout}[1]{}
\newcommand{\cZ}{{\cal Z}}

\newcommand{\dx}{\mathop{}\! dx}
\newcommand{\dy}{\mathop{}\! dy}

\begin{document}

\author{Xaver Kriechbaum\footnote{Department of Mathematics, Weizmann Institute of Science, Rehovot 76100, Israel; xaver.kriechbaum@weizmann.ac.il} \and Lenya Ryzhik\footnote{Department of Mathematics, Stanford University, Stanford, CA 94305, USA; ryzhik@stanford.edu}
 \and Ofer Zeitouni\footnote{Department of Mathematics, Weizmann Institute of Science, Rehovot 76100, Israel; ofer.zeitouni@weizmann.ac.il}}
\title{Voting models and tightness for a family of recursion equations}

\maketitle

\begin{abstract} 
We consider recursion equations of the form
$u_{n+1}(x)=Q[u_n](x),~n\ge 1,~x\in\Rm$,
with a non-local operator 
$Q[u](x)= g( u\ast q)$, where $g$ is a polynomial, satisfying
$g(0)=0$, $g(1)=1$, $g((0,1)) \subseteq (0,1)$, and $q$ is a (compactly supported) probability density with
$\ast$ denoting convolution. Motivated by a line of works for
nonlinear PDEs initiated by Etheridge, Freeman and Penington (2017), 
we show that for general $g$, a probabilistic 
model based on branching random walk
can be given to the solution of the recursion,  while in case $g$ is also strictly monotone, a probabilistic threshold-based model 
can be given. In the latter case, we provide a conditional tightness result. We analyze in detail the 
bistable case
and prove for it convergence of the solution shifted  around a linear in $n$
centering.
\end{abstract}

\section{Introduction}
We consider in this paper certain recursion equations that are discrete-time 
analogs of the (nonlinear) PDE
\begin{equation}
  \label{eq-1611a}
  \partial_t u(t,x)=\frac12 \partial_{xx} u(t,x)+f(u),\quad
  x\in \Rm,~ t\in \Rm_+.
\end{equation}
Here, $f$ is (typically) a polynomial satisfying $f(0)=f(1)=0$, and $u(t,x)$
is assumed to satisfy the boundary conditions 
\begin{equation}
  \label{eq-1611b}
  \lim_{x\to -\infty} u(t,x)=1, \quad \lim_{x\to\infty} u(t,x)=0.
\end{equation}

An important special case is $f(u)=u-u^2$, when \eqref{eq-1611a} is the
so called Fisher--Kolmogorov-Petrovskii-Piskunov (FKPP) equation
\cite{Fisher,KPP}. Then, (\ref{eq-1611a})-(\ref{eq-1611b}) 
admits traveling wave solutions of the 
form~$u(t,x)=w(x-vt)$ for all $v\ge v_*=\sqrt{2}$,  and the solution 
to \eqref{eq-1611a} with an initial condition that is compactly supported on the right, after proper centering,
converges to the traveling wave \cite{KPP} moving with the minimal
speed $v=v_*$. In a celebrated work, Bramson 
\cite{Bramson83} computed the centering.
An important observation, often attributed to McKean 
\cite{Mckean} but going back at least 
to Skorohod~\cite{Skorohod}, gives a representation of the solution of \eqref{eq-1611a} 
with step initial condition $u(0,x)=\one({x<0})$,
in terms of a branching Brownian motion. It is defined as follows: start with a particle at the origin that performs a Brownian motion. At an independent, exponentially distributed time $\tau$, the particle splits in two, and each particle starts afresh and independently, from its current location, the same process. With $N_t$ the number of particles at time $t$, and with $(X_t^i)_{i=1,\ldots, N_t}$ denoting their positions and $M_t=\max_{i} X_t^i$, we have that $u(t,x)=\Pm(M_t\geq x)$.
In particular, that representation
is at the heart of Bramson's computation of the centering term.

There is an analogous story for discrete recursions. Namely, consider the 
recursion
\be\label{sep1402}
u_{n+1}(x)=Q[u_n](x),~~n\ge 1,~x\in\Rm,
\ee
with a non-local operator 
\be\label{eq-1611c}
Q[u](x)= g( u\ast q),
\ee
where $\ast$ denotes convolution, $g(x)=f(x)-x$,
and $q$ is a probability measure. (We refer to \cite{AlB} for a general discussion of such recursions.) In the particular
case $g(x)=x^2$, one has a probabilistic interpretation of the solution
in terms of the law of the
maximum displacement of a branching random walk (BRW)
with binary branching and increment law $q$. For such BRWs, convergence of the
law of the centered maximum, evaluation of the centering, and identification
of the limit, were obtained by A\"{i}d\'{e}kon
\cite{Aidekon}, see also \cite{BDZ}, after some initial results on tightness were described in \cite{ABR} and \cite{BramZeiRec}.

Returning to the PDE setup, the convergence to a traveling wave extends to 
a family of ``KPP-like'' nonlinearities, which in particular
do not possess any zero in the interval $(0,1)$. In case such zeroes exist, 
some partial results are contained in \cite{FifeMcL} (for the so called bistable case)
and \cite{gilettimatano}; in general, convergence to a traveling
wave is replaced by the notion of existence of ``terraces'', of 
increasing width and connected by travelling waves.

Still in the context of PDEs,
the probabilistic representation of Skorohod and McKean extends readily 
to the situation where 
\begin{equation}
  \label{eq-1711a}
  f(u)=1-\sum_k p_k (1-u)^k -u
\end{equation}
with $p_k\geq 0$ and 
$\sum p_k=1$,
by modifying the branching mechanism from binary to random with law $p_k$.
This can be further extended 
to a limited class of nonlinearities $f$ of that type, 
see \cite{Watanabe,Japanese}.  

A major breakthrough concerning probabilistic representations for the solutions to
\eqref{eq-1611a} came with the work \cite{EFP}. Motivated by the Allen-Cahn
equation, it deals with the nonlinearity 
\[
f(u)=u(1-u)(2u-1),
\]
and proposed
a probabilistic representation based on BBM with ternary branching followed
by a ``voting rule'' that propagates the locations of the particles at
time $t$ through the genealogical tree 
to a random variable, whose law represents the solution.
That this representation applies to arbitrary polynomial $f$ was 
observed shortly after in \cite{odowd} and \cite{an22voting}.

Concerning the discrete setup, for  nonlinearities of the form 
\eqref{eq-1711a}, a certain steepness comparison present in the continuous 
setup does not transfer to the discrete case unless the density $q$
is log-concave; see \cite{bachmann}. For more general densities of compact support, a clever
probabilistic argument that yields tightness was
presented in \cite{DH91}, 
while an analytic argument,
based on the recursions \eqref{eq-1711a} and applying to a wide class of positive $f$ under mild assumptions on $q$,
was presented in \cite{BramZeiRec}.

Our goal in this paper is to study
the discrete recursions \eqref{sep1402} with polynomial functions $g$,
and develop for them a probabilistic 
representation similar to that studied in 
\cite{EFP,odowd,an22voting}. As in \cite{an22voting}, we distinguish between random
threshold
models and random outcome models, and show in Propositions~\ref{Theo:ReprOutc}
and~\ref{Theo:NonLinThresh}
that to any polynomial
$g$ with $g(0)=0$, $g(1)=1$, $g((0,1))\subseteq (0,1)$ one can find a random outcome
model which represents the solution to \eqref{sep1402}, 
while a random threshold model can be found only if $g$ is, in addition,
monotone (note that $f(x)=g(x)-x$ is not required to be monotone). 
In the latter case, we use the probabilistic representation and a modification of the Dekking-Host argument to prove in Theorem 
\ref{Theo:Clustering}  the existence of terraces, interpreted as conditional 
tightness statements; 
we also analyze in some details the case of 
binary-ternary branching with threshold voting, 
see  Section \ref{subsec-bt}. We chose to do so because of the very clear 
probabilistic interpretation of the voting rule in that particular model 
(see the min-max~$\widetilde{M}_n$ in
\eqref{eq-1711d}),
and because standard techniques for handling the maximum of BRW do not seem
to work for handling the min-max $\widetilde{M}_n$.
Section 
\ref{Sec:TightnessAnMeans} is devoted to an analytical study of the bistable
case (where $f(x)=g(x)-x$ possesses a single zero in $(0,1)$); convergence to a travelling wave (with linear in $n$ centering) is proved.

\subsection{Notation and setup}

Throughout, $q$ denotes a probability measure on $\Rm$ which we assume
to possess 
a density $q(\cdot)$ with respect to the Lebesgue measure.
We further  assume that the density~$q(x)$ is continuous and has compact 
support, and fix $C_q>0$ such that
\be\label{aug3025}
q\in C_c(\Rm),~~\hbox{supp}(q)\subseteq[-C_q,C_q].
\ee


The density $q$ will serve as
the jump density of the increments of a branching random walk,
with offspring law $\{p_d\}$; that is,
$p_d$ denotes the probability for a parent to have $d$ children. The 
resulting rooted 
Galton--Watson tree up to generation $n$ is denoted $\cT_n$, with $o$ denoting the root; explicitly,
$\cT_n$ is a random tree with vertex set still denoted
$\cT_n$, and edge set $E_n$. For a vertex $v\in \cT_n$, we denote by $|v|$ 
its (tree) distance from the root,
and we let $D_n=\{v\in \cT_n: |v|=n\}$,
which corresponds to the collection of particles at time $n$.
We denote by $S_v^z$, $v\in \cT_n$, a BRW that starts at $z\in \Rm$.
We note that under our assumptions on $q$, 
for each $x,z\in\Rm$ and any vertex $v\neq o\in\cT_n$ we have
\be\label{aug2412}
\Pm\big[S_v^z=x\big]=0.
\ee

Given a collection of numbers $x_1,\dots,x_n$ and $k\le n$, we denote by $x_{(k)}$ the $k$-th largest element in that collection,
so that $x_{(1)}\le x_{(2)}\le \dots\le x_{(n)}$. 

{\bf Acknowledgement.} The work of XK and OZ 
was supported by Israel Science Foundation grant number 421/20.
LR was supported by NSF grants DMS-1910023 and DMS-2205497 and by~ONR grant N00014-22-1-2174.  
We thank Alison Etheridge and Jean-Michel Roquejoffre for useful discussions.

\section{Recursions as voting models for branching random walks}

In this section, we first define the discrete analogs to the random outcome and random threshold voting models, as defined in \cite{an22voting}. 
After this, we discuss which nonlinearities can be achieved in the recursions associated to these models. One surprising difference to the continuous model is that in the discrete 
case the random outcome model is more general in the sense that there are nonlinearities we can describe using it, which can
not be described with the random threshold model. We should note that  the probabilistic side of the analysis in this paper will only work for the random threshold model.

\subsection{Voting models and recursive equations}

\subsubsection{Random threshold models}\label{sec:ran-thresh}

A random threshold voting model on the Galton-Watson tree ${\cal T}_n$ of the branching random walk $S_v^0$ with $v\in \cT_n$ 
is defined as follows. First, at the final time $n$ we assign the values $\vphi_n(v)=S_v^0$ for all vertices~$v\in D_n$.
Next, at each vertex $v$ of the tree $\cT_n$ with $|v|<n$, let 
\be\label{aug2408}
d(v)=|D_1(v)|,
\ee
be the number of children of the vertex $v$. 
Then, we choose a number $L_v\in[1,2,\dots,d(v)]$, with the probabilities 
\be\label{aug2402}
\mathbb{P}\big[L_v = k | d(v)=d\big] = \zeta_{k,d}.
\ee
Here,  $\zeta_{k,d}\in[0,1]$ are assigned, so that
\be\label{aug2404}
\sum_{k=1}^d \zeta_{k,d} = 1,~~\hbox{for all $d$}.
\ee
We can now propagate the values $\vphi_n(v)$ up the genealogical tree $\cT_n$ recursively, 
by assigning to a given vertex $v$ with $|v|<n$ the value $\vphi_n(v)$ that is the $L_v$-th largest of the values of $\vphi_n(w)$, where $w$ are
all the children
of $v$. That is, if we order  $w_j\in D_1(v)$, $j=1,\dots,d(v)$, 
according to the increasing order of $\varphi(w_j)$, then 
\be\label{aug2406}
\vphi_n(v)=\vphi_n(w_{L_v}),
\ee
Finally we set 
\be\label{aug2410}
M_n := \varphi_n(o).
\ee

Equivalently, we can consider a  voting process, for a BRW $S_v^x$ that starts originally at a position $x\in\Rm$. The particles 
$v$ in generation 
$n$ vote $1$ if and only if $S_v^x\ge 0$ and a particle in a generation $k<n$ votes~$1$ if less than $L_v$ of its children voted $0$. 
From the construction, it is immediate to see that a particle~$v\in\cT_n$ 
votes $0$ iff $\varphi_n(v)<0$. Thus, using $V_n^x(o)$ to denote the vote at the origin when the BRW starts at $x$, we have for $n\ge 1$
\[
\mathbb{P}[M_n\le x] = \mathbb{P}[M_n<x] = \mathbb{P}[V_n^{-x}(o) = 0].
\]
In the first step above, we used (\ref{aug2412}) that gives
\[
\mathbb{P}[M_n = x] \le \mathbb{P}\Big[\bigcup_{v  : |v| = n} \{S_v = x \}\Big] = 0.
\]

It is straightforward to use the definition of $M_n$ and the independence of the increments of BRW, 
to deduce that the distribution function 
\be
F_{n}(x)=\Pm[M_n\leq x]
\ee
satisfies the renewal equation
\begin{equation}
\begin{aligned}
&F_{{n+1}} = g(q\ast F_{n}), \\ 
&F_{0} = \one(x\ge 0), 
\end{aligned} \label{eq:RecEqThresh}
\end{equation}
with
\begin{equation}\label{eq:gThresh}
g(x) = \sum_{d=1}^{d_0} \sum_{k=1}^d p_d\zeta_{k,d} \sum_{l=k}^d \binom{d}{l}x^l (1-x)^{d-l}. 
\end{equation}
We will write 
\be\label{aug2802}
f(u)=g(u)-u,
\ee
and call $f$ the nonlinearity and $g$ the recursion polynomial associated to the random threshold model.
Note that up to the prefactor $\beta$, the function $f$ coincides with equation (3.35) in \cite{an22voting}. We remark that in the special case when
one chooses $L_v = |D_1(v)|$ deterministically, 
the value $M_n$ is the maximum of the underlying BRW in generation $n$. In this sense, random threshold models are a generalization of the study of the maximum of BRWs.

\subsubsection{Random outcome models}

In contrast, a random outcome voting model is defined as follows. Let $\cT_n$ be the genealogical tree of a BRW that originally starts at a position $x\in\Rm$.
We  fix the probabilities $\alpha_{k,d}\in[0,1]$, defined for $d\ge 1$ and~$0\le k\le d$, such that 
\be
\alpha_{0,0}=0,~~ \alpha_{d,d}=1.
\ee
The voting on $\cT_n$ is done as follows. For a final generation particle $v$, such that $|v|=n$, we set  
\be
V_n^x(v) := \one(S_v^x\ge 0).
\ee
For a vertex $v$ with $|v|<n$, such that $k$ out of its $d$ children voted one, 
we let $V_n^x(v)$ be a random variable with 
\be
\mathbb{P}[V_n^x(v) = 1] = \alpha_{k,d}.
\ee
Then, the function $u_n(x) = \mathbb{P}[V_n^{x}(o) = 1]$ satisfies the recursion equation
\begin{equation}
\begin{aligned}
&u_{n+1}= g(\widehat{q}\ast u_n),\\
&u_0 = \one(x\ge 0),
\end{aligned} \label{eq:RecOutc}
\end{equation}
where $\widehat{q}(x) = q(-x)$ and 
\begin{equation}
g(x) = \sum_{d=1}^{d_0}\sum_{k=1}^d p_d\alpha_{k,d} \binom{d}{k}x^k(1-x)^{d-k}.\label{eq:gOutc}
\end{equation}
We call $g$ the recursion polynomial associated to the random outcome model. We will see in Section \ref{Sec:WhatNon} that random outcome models are a further generalisation of random threshold models.

\subsection{Background on the Bernstein polynomials}

The recursion polynomials coming from a voting scheme are convenient to represent in terms of the Bernstein polynomials
\be
b_{k,d}(x) := \binom{d}{k} x^k(1-x)^{d-k}.
\ee
We use the convention $b_{k,d} \equiv 0$ if $k\not\in\{0,\dots, d\}$. In this section, we recall several useful properties of the Bernstein polynomials.

First, we note that the Bernstein polynomials of degree $d$ form a basis of the space $\R^{\le d}[x]$ of the polynomials of degree lesser equal $d$. 
For  a polynomial $p(x)$ we denote the coefficients with regard to the Bernstein polynomials of degree $d\ge \mathrm{deg}(p)$ by $\beta_{k,d}(p)$:
\be
p(x)=\sum_{k=0}^d\beta_{k,d}(p)b_{k,d}(x).
\ee
Additionally, as $b_{k,d}(0) = 0$ for $k\neq 0$ and $b_{k,d}(1) = 0$ for $k\neq d$,  it follows 
that $(b_{1,d},\dots, b_{d-1,d})$ form a basis of  the sub-space $\{p\in \R^{\le d}[x]: p(0) = p(1) = 0\}$.  

The Bernstein polynomials satisfy the following elementary  algebraic identities. First,  for all $d\in \N$ we have
\begin{align}
\sum_{k=0}^d b_{k,d}(x)\equiv 1.
\label{Def:BernSum}
\end{align}
Second,  for all $d\in\N$,  we have
\be\label{eq:BernDeriv}
\bal
b_{k,d}'(x) &= d(b_{k-1,d-1}(x)-b_{k,d-1}(x)),~~\hbox{for $1\le k\le d-1$},\\
b_{0,d}'(x)&=-db_{0,d-1},~~b_{d,d}'=db_{d-1,d-1},
\enbal
\ee
and
\be \label{eq:BernDegEl}
b_{k,d-1}(x) = \frac{d-k}{d}b_{k,d}(x)+\frac{k+1}{d}b_{k+1,d}(x),~~\hbox{ for $0\le k\le d-1$.}
\ee

Next, we recall a way to compute the coefficients of a polynomial $p(x)$ in the Bernstein basis of degree~$d+1$ from the coefficients in degree $d\ge \mathrm{deg}(p)$, 
compare to equation (12) in \cite{BernstPoly}:
\begin{equation}
\beta_{k,d+1}(p) =
\begin{cases}
\beta_{0,d}(p),&\quad \text{for}\ k = 0,\\
\frac{k}{d+1}\beta_{k-1,d}(p)+\frac{d+1-k}{d+1}\beta_{k,d},&\quad \text{for}\ 1\le k\le d,\\
\beta_{d,d}(p),&\quad \text{for}\ k = d+1.\\
\end{cases} \label{eq:ElevateCoeff}
\end{equation}

Finally, we cite two results about getting bounds on $\beta_{k,d}(p)$ from bounds on $p$.
\begin{prop}[Theorem 2 in \cite{BernstPoly}]\label{Lem:BernstRepr}
Given a polynomial $p(x)$, there exists $d\ge\mathrm{deg}(p)$ such that the Bernstein coefficients $\beta_{k,d}(p)$ satisfy $0\le\beta_{k,d}\le 1$
for all $0\le k\le d$ if and only if either (i) $p(x)\equiv 0$, or (ii) $p(x)\equiv 1$, or (iii) $0\le p(0),p(1)\le 1$ and~$0<p(x)<1$ for all $x\in(0,1)$.
\end{prop}
\begin{prop}\label{Lem:BernCoeffPos}
Let $p(x)$ be a polynomial such that $p(x)>0$ for all $x\in (0,1)$. Then there is~$d_0\ge\deg(p)$ such that for all $d\ge d_0$ and all $k\in\{0,\dots, d\}$ we have $\beta_{k,d}(p)\ge 0$. 
\end{prop}
While Proposition \ref{Lem:BernCoeffPos} can't be found verbatim in \cite{BernstPoly} it can be easily recovered from the proof of their Theorem 4.

\subsection{Achievable recursions} \label{Sec:WhatNon}

In this section we  explain which recursions can be represented with a random threshold  or a random outcome model. 
One notable difference to the corresponding Theorems 3.2 and 3.3 from \cite{an22voting} for voting models for a branching Brownian motion 
is that the random threshold model can represent (strictly) less recursions than the random threshold model. The first
result characterizes the polynomials that can be represented via a random outcome model.  This result is very similar to
Theorem 3.2 of  \cite{an22voting}. 
\begin{prop}\label{Theo:ReprOutc} Let $g(x)$ be a polynomial. The following are equivalent:\\
(i) there is a random outcome model with recursion polynomial $g(x)$,\\
(ii) there is $d\ge \deg(g)$ and $\alpha_{k,d}$, $0\le k\le d$ such that $\alpha_{0,d}=0$, $\alpha_{d,d}=1$, $0\le\alpha_{k,d}\le 1$ for
all~$1\le k\le d-1$, and
\be
g(x)=\sum_{k=0}^d\alpha_{k,d}b_{k,d}(x),
\ee
(iii) $g(0)=0$, $g(1)=1$ and $0<g(x)<1$ for all $x\in(0,1)$.
\end{prop}
\begin{proof}
We denote the set in (i)  by $V_1$, the one in (ii)  by $V_2$ and the one in (iii) by $V_3$. The fact that~$V_2 = V_3$
 is an immediate consequence of Proposition~\ref{Lem:BernstRepr} and the observation that if $\alpha_{k,d}$ are the Bernstein coefficients 
 of $g(x)$, then $g(0) = \alpha_{0,d}$, $g(1) = \alpha_{d,d}$.

The inclusion $V_2\subseteq V_1$ follows from \eqref{eq:gOutc} by considering a BRW with $p_k =1$ if $k=d$ and $p_k=0$ otherwise.

Finally, the inclusion $V_1\subseteq V_3$ follows from \eqref{eq:gOutc}, the requirement that $\alpha_{0,d} = 0$, $\alpha_{d,d} = 1$,
$0\le\alpha_{k,d}\le 1$ for $1\le k\le d-1$ 
and~(\ref{Def:BernSum}).
\end{proof}

The next result characterizes the polynomials that can be represented by a random threshold model. The result is different from the continuous case~\cite{an22voting}.
The reason is that in the continuous case such representations may require a very fast exponential clock, which we do not have available for BRW.  
\begin{prop}\label{Theo:NonLinThresh}
Let $g(x)$ be a polynomial. The following are equivalent:\\
(i) there is a random threshold model with recursion polynomial $g(x)$,\\
(ii) there is $d\ge \deg(g)$ and $\alpha_{k,d}$, $0\le k\le d$ 
such that 
\be\label{aug2416}
\hbox{$\alpha_{0,d}=0$, $\alpha_{d,d}=1$ and $0\le\alpha_{k-1,d}\le\alpha_{k,d}\le 1$,  for
all~$1\le k\le d-1$,}
\ee
and
\be
g(x)=\sum_{k=0}^d\alpha_{k,d}b_{k,d}(x),
\ee
(iii) $g(0)=0$, $g(1)=1$, and both $0<g(x)<1$ and $g'(x)>0$ for all $x\in(0,1)$.
\end{prop}
\begin{proof}
We denote the set in  (i) by $W_1$, the one in (ii)  by $W_2$ and the one in (iii) by $W_3$.

We start by proving that $W_2 \subseteq W_3$. Given $g\in W_2$, we know that 
\be\label{aug2502}
g(0)=\alpha_{0,d}=0,~~g(1)=\alpha_{d,d}=1.
\ee
Next, we use  \eqref{eq:BernDeriv} and (\ref{aug2502}), to obtain
\begin{equation} \label{eq:gderiv}
\bal
g'(x) &= \sum_{k=1}^d\alpha_{k,d}b_{k,d}'(x)=d\sum_{k=1}^{d-1} \alpha_{k,d}(b_{k-1,d-1}-b_{k,d-1})+db_{d-1,d-1}\\
&=d\sum_{k=0}^{d-1}\alpha_{k+1,d}b_{k,d-1}(x)-\sum_{k=1}^{d-1}\alpha_{k,d}b_{k,d-1}
= d\sum_{k=0}^{d-1} (\alpha_{k+1,d}-\alpha_{k,d})b_{k,d-1}(x). 
\enbal
\end{equation}
We see from (\ref{aug2416}) that each term in the last sum above is non-negative  and there has to be least one~$k$ 
such that~$\alpha_{k+1,d}-\alpha_{k,d}>0$. As $b_{k,d-1}(x)>0$ for all $x\in (0,1)$, it follows that $g'(x)>0$ for all $x\in (0,1)$.

Next we prove that $W_3\subseteq W_2$. Given $g\in W_3$, by Proposition \ref{Lem:BernstRepr} there is 
$d_0$ such that 
\[
0\le\beta_{k,d_0}(g)\le 1,~~\hbox{for all $0\le k\le d_0$},
\]
while we also have 
\be\label{aug2504}
\beta_{0,d}(g) = g(0) = 0,~~\beta_{d,d}(g) = g(1)  = 1.
\ee 
Using \eqref{eq:ElevateCoeff} yields 
that for all~$d\ge d_0$ we also have 
\be\label{aug2506}
0\le\beta_{k,d}(g)\le 1,~~\hbox{for all $0\le k\le d$},~~\beta_{0,d}(g) = 0,~~\beta_{d,d}(g) = 1.
\ee
Applying Proposition \ref{Lem:BernCoeffPos}, we deduce, in addition, that there is $d_1$ such that for all $d\ge d_1$ and $0\le k\le d$ 
we have~$\beta_{k,d}(g')\ge 0$. Then, for any  $d\ge \max\{d_0,d_1+1\}$ we have, using (\ref{eq:BernDeriv}) and (\ref{aug2506})
\begin{align*}
g'(x) &= \sum_{k=0}^d \beta_{k,d}(g) b_{k,d}'(x) =\sum_{k=1}^d \beta_{k,d}(g) b_{k,d}'(x) 
=d\sum_{k=1}^{d-1}\beta_{k,d}(g)[b_{k-1,d-1}(x)-b_{k,d-1}(x)] +db_{d-1,d-1}(x)\\
&=d\sum_{k=0}^{d-1}\beta_{k+1,d}(g)b_{k,d-1}(x)-d\sum_{k=0}^{d-1}\beta_{k,d}(g)b_{k,d-1}(x)
= d\sum_{k=0}^{d-1} (\beta_{k+1,d}(g)-\beta_{k,d}(g))b_{k,d-1}(x).
\end{align*}
This implies that 
\be
d(\beta_{k+1,d}(g)-\beta_{k,d}(g)) = \beta_{k,d-1}(g')\ge 0,~~\hbox{for all $0\le k\le d-1$.}
\ee 
Here, the last step used $d-1\ge d_1$. Thus, we have 
\be
\beta_{0,d}(g)\le\beta_{1,d}(g)\le \dots \le \beta_{d,d}(g),
\ee 
which is the second condition in (\ref{aug2416}). Combined with (\ref{aug2506}) that holds because $d\ge d_0$, 
we see that the first condition in (\ref{aug2416}) also holds, and $g\in W_2$.

Next we show that $W_2\subseteq W_1$. Take $g\in W_2$ and write it as 
\be
g(x)= \sum_{k=0}^d \alpha_{k,d} b_{k,d}(x),
\ee
with $\alpha_{k,d}$ as in (\ref{aug2416}). Consider a BRW with branching into $d$ children and a random threshold voting model with 
\be\label{aug2508}
\zeta_{0,d}=0,~\zeta_{k,d} = \alpha_{k,d}-\alpha_{k-1,d},~~\hbox{ for $1\le k\le d$.}
\ee
Since 
\be
\sum_{k=1}^d \zeta_{k,d} = \alpha_{d,d}-\alpha_{0,d} = 1,
\ee
this does define a random threshold model. By \eqref{eq:gThresh},  the associated recursion polynomial is
\begin{align*}
\widetilde{g}(x) = \sum_{k=1}^{d} \zeta_{k,d}\sum_{l=k}^db_{l,d}(x) = \sum_{k=1}^d \left(\sum_{l=1}^k \zeta_{l,d}\right)b_{k,d}(x) 
= \sum_{k=1}^d \left(\alpha_{k,d}-\alpha_{0,d}\right)b_{k,d}(x) =\sum_{k=0}^d\alpha_{k,d}b_{k,d}(x)= g(x).
\end{align*}
Here, the last step used  $\alpha_{0,d} = 0$. 

Finally we prove that $W_1\subseteq W_3$. Given $g\in W_1$, it has the form (\ref{eq:gThresh}):
\be\label{aug2510}
g(x) = \sum_{d=1}^{d_0}\sum_{k=1}^d p_d\zeta_{k,d} \sum_{l=k}^d b_{l,d}(x),
\ee
so that 
\be\label{aug2514}
g(0) = \sum_{d=1}^{d_0}\sum_{k=1}^d p_d\zeta_{k,d} \sum_{l=k}^d b_{l,d}(0) =0,
\ee
since $b_{k,d}(0) = 0$ for $k\neq 0$. Furthermore, we have
\be\label{aug2516}
g(1) = \sum_{d=1}^{d_0}\sum_{k=1}^dp_d\zeta_{k,d}\sum_{l=k}^db_{l,d}(1) = \sum_{d=1}^{d_0}p_d\sum_{k=1}^d \zeta_{k,d} = 1,
\ee
since $b_{d,d}(1) = 1$, $b_{k,d}(1) = 0$ for $k\neq d$,  and
\be\label{aug2512}
\sum_{k=1}^d\zeta_{k,d} = \sum_{d=1}^{d_0} p_d = 1.
\ee
We also note that (\ref{Def:BernSum}) and  (\ref{aug2512}) imply that
\be
0<g(x)< \sum_{d=1}^{d_0}\sum_{k=1}^d p_d\zeta_{k,d} \sum_{l=0}^d b_{l,d}(x)=1,~~\hbox{ for all $0< x< 1$.}
\ee
Finally, using \eqref{eq:BernDeriv} yields that
\be\label{aug2518}
g'(x) = \sum_{d=1}^{d_0} \sum_{k=1}^dp_d\zeta_{k,d}\Big[\sum_{l=k}^{d-1} d(b_{l-1,d-1}-b_{l,d-1}) +db_{d-1,d-1}\Big] = \sum_{d=1}^{d_0} \sum_{k=1}^d dp_d\zeta_{k,d} b_{k-1,d-1}(x)>0,
\ee
for all $x\in (0,1)$. Combining (\ref{aug2514}), (\ref{aug2516}) and (\ref{aug2518}), we conclude that $g\in W_3$.  
\end{proof}

\section{Clustering with probabilistic means}


In this section, we consider a random threshold voting model, as described in Section ~\ref{sec:ran-thresh}, 
and the corresponding ``result of the voting" $M_n$, defined by (\ref{aug2410}). Our goal is to
prove that there are intervals~$I_{j,n}$, $j=1,\dots,N_I$, such that $\bigcup_{j=1}^{N_I}I_{j,n} = \R$ and $M_n$ conditioned to stay inside $I_{j,n}$ is tight around its median for all $j\in\N$. 

Let $f(u)$ be the nonlinearity coming from a random threshold voting model, as defined by~\eqref{eq:gThresh}-(\ref{aug2802}), 
and $N_f$ be the number of zeroes of $f(u)$ inside the open interval $(0,1)$. 
We will see that the cumulative distribution $F_{M_n}$ of~$M_n$ is composed out of at most $N_f+1$ clusters. 
However,  we cannot tell whether some of these clusters coincide. 
In particular we cannot deduce how many ``terraces''~$F_{M_n}$ has, but only give an upper bound on their number. 
If there is just one cluster, then there is a  sequence~$m_n$ such that~$(M_n-m_n)_{n\in\N}$ is tight.
In particular, we reprove (in the case of compact support for the increments)
the result from \cite{BramZeiRec} that for~$f$ with no zeroes inside~$(0,1)$, 
the sequence $(M_n-\mathrm{med}(M_n))_{n\in\N}$ is tight. 

For $n\in\Nm$ fixed, given an interval $I$, we let $M_{n;I}$ be a random variable with 
\be
\mathbb{P}[M_{n;I}\in \cdot] := \mathbb{P}[M_n\in \cdot | M_n\in I].
\ee
Here is the main result of this section. 
\begin{thm}\label{Theo:Clustering}
Let $(p_d,\zeta_{k,d})_{d\le d_0, k\in\{1,\dots, d\}}$ be a random threshold model such that the associated nonlinearity $f \not\equiv 0$. 
Let $0 = \alpha_0<\alpha_1<\dots<\alpha_{N_f+1}  =1$ be the zeroes of $f$, and $q_{n,s}$, $s\in\{0,\dots, N_f+1\}$, be the $\alpha_s$-quantile of $M_n$,
and set $I_{s,n} := [q_{s-1,n}, q_{s,n}]$. Then, for all $s\in \{1,\dots, N_f+1\}$ the sequence $(M_{n;I_{s,n}}-\mathrm{med}(M_{n;I_{s,n}}))_{n\in\N}$ is tight.
\end{thm}

The assumption $f\not\equiv 0$ is necessary. This can be seen by considering the voting model for any BRW with the probabilities 
$\zeta_{k,d} = {1}/{d}$ for all $d\le d_0$, $k\in\{1,\dots, d\}$. 
This corresponds to the parent using the value of one of its children uniformly at random, which means that the particle in the $n$-generation whose value 
is propagated to the top is chosen uniformly and random. Then, we have $M_n \stackrel{d}{=} S_n$, and $M_n$ is not tight but has distribution function which 
spreads as $\sqrt{n}$.
 
Heuristically the proof of Theorem \ref{Theo:Clustering} shows that $M_{n;I}$ has a recursive structure similar to \eqref{eq:RecEqThresh} but governed by the 
recursion associated to the nonlinearity $f_{|_I}$ rescaled to be a function with domain $[0,1]$, defined in (\ref{aug2910}) below. 
Since this nonlinearity doesn't have a zero in $(0,1)$ we get tightness similar to the soft argument for tightness of the maximum of BRWs with bounded 
increments given in \cite{DH91}.

\begin{figure}
\centering
\begin{subfigure}{0.3\textwidth}
\includegraphics[width = \textwidth]{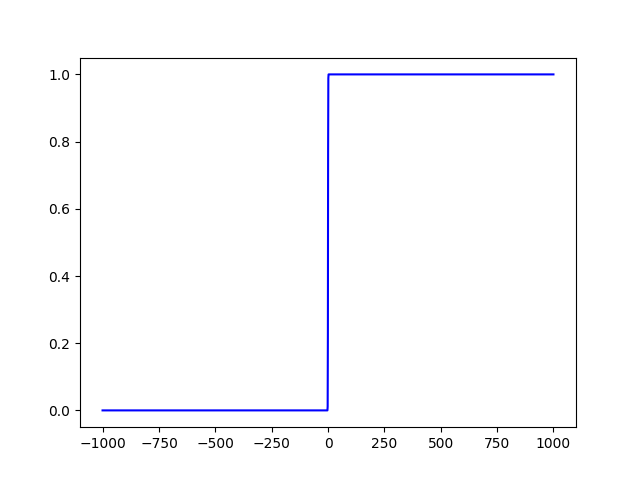}
\caption{Here $\zeta_{3,4} = 1/2$, $\zeta_{2,4} = 1/2$. \\\strut}
\label{fig:Tight}
\end{subfigure}
\hfill
\begin{subfigure}{0.3\textwidth}
\includegraphics[width = \textwidth]{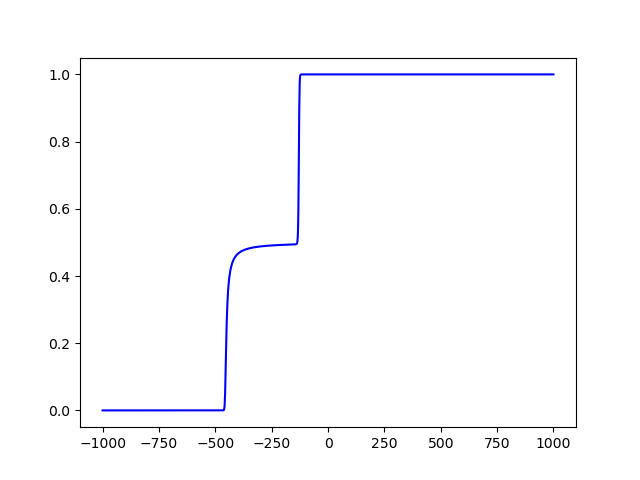}
\caption{Here $\zeta_{4,4} = 3/16$, $\zeta_{3,4} = 19/48$, $\zeta_{2,4} = 5/48$ and $\zeta_{1,4} = 5/16$.}
\label{fig:OneSideTight}
\end{subfigure}
\hfill
\begin{subfigure}{0.3\textwidth}
\includegraphics[width = \textwidth]{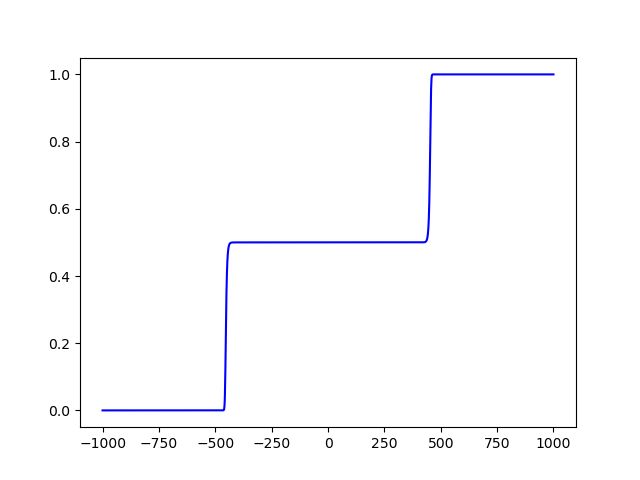}
\caption{Here $\zeta_{4,4} = 5/16$, $\zeta_{3,4} = 3/16$, $\zeta_{2,4} = 3/16$ and $\zeta_{1,4} = 5/16$.}
\label{fig:NonTight}
\end{subfigure}
\caption{The plots show distribution functions of $M_{1000}$ for $q = 1/4\delta_{-1}+1/4\delta_{1}+1/2\delta_0$, $p_4 = 1$ and various $\zeta_{k,4}$, which are specified above.}
\label{fig:1}
\end{figure}

To illustrate Theorem \ref{Theo:Clustering}, let us consider  the case when $f$ has a single zero $\alpha$. We denote by $q_n$ the~$\alpha$-quantile of $M_n$.
The distribution of $M_n$ has three possible archetypes. In the three examples below, we consider BRW with branching into four children,
that is, $p_4=1$:
\begin{enumerate}[(i)]
\item The distribution has two tight clusters, which are $O_n(1)$ far away from $q_n$, in particular $(M_n-q_n)$ is tight. An example  of this behavior comes from
$\zeta_{3,4} = 1/2$, $\zeta_{2,4} = 1/2$.
The distribution function of $M_{1000}$ can be seen in Figure \ref{fig:1} \subref{fig:Tight}.
\item The distribution has two tight clusters, one of which is at distance $O_n(1)$ to $q_n$ and one which  moves away from $q_n$. 
An example for this  is $\zeta_{4,4} = 3/16$, $\zeta_{3,4} = 19/48$, $\zeta_{2,4} = 5/48$ and $\zeta_{1,4} = 5/16$. The distribution function of $M_{1000}$ can be seen in Figure \ref{fig:1} \subref{fig:OneSideTight}. We note that such examples cannot be generated with a symmetric voting rule, since in the symmetric situation $M_{n;I_1} \stackrel{d}{=} -M_{n;I_2} \stackrel{d}{=} -|M_n|$. 
\item The distribution has two tight clusters, both of which are further than distance $O_n(1)$ away from~$q_n$. An example for this is  
$\zeta_{4,4} = 5/16$, $\zeta_{3,4} = 3/16$, $\zeta_{2,4} = 3/16$ and $\zeta_{1,4} = 5/16$, the distribution function of $M_{1000}$ can be seen in Figure \ref{fig:1} \subref{fig:NonTight}.
\end{enumerate}

Before we state the main ingredient in the proof of Theorem \ref{Theo:Clustering}, let us introduce some notation. 
For~$s\in \{1,\dots, N_f+1\}$, we define a stretched version of the restriction of $f(x)$ to $I_{s,n}$:
\be\label{aug2910}
\widetilde{f}_s(x) := f((\alpha_s-\alpha_{s-1})x+\alpha_{s-1}),
\ee
so that
\be\label{aug2904}
\widetilde f_s(0)=\widetilde f_s(1)= 0,~~\widetilde f(x)\neq 0,~~\hbox{ for all $x\in(0,1)$.}
\ee
We also define a piece-wise linear function
\be\label{aug2914}
\psi_{s,n} (x)= (q_{s,n}-q_{s-1,n})\one(x\le q_{s-1,n}) +(q_{s,n}-x)\one({x\in [q_{s-1,n},q_{s,n}]}) .
\ee
Finally, for each $n\in\N$, we denote by $M_{n,k}$,  $k\in\N$, a family of i.i.d.~random variables such that  $M_{n,1}\stackrel{d}{=}M_n$,
and let $(M_{n,k;I})_{k\in\N}$ be i.i.d.\@ with $M_{n,1;I} \stackrel{d}{=} M_{n;I}$. In addition, given any $D\ge 1$, we let $M_{n,(k);I}^{(D)}$ 
be the $k$-th largest element of $(M_{n,k;I})_{k\in(1,\dots,D)}$, so that
\be\label{aug2912}
M_{n,(1);I}^{(D)}\le M_{n,(2);I}^{(D)}\le\dots\le M_{n,(D);I}^{(D)}.
\ee
\begin{lem}\label{Lem:MainIngrClust}
We have, for all $D\ge d_0$
\be\label{aug2916}
-C_q\le \mathbb{E}[\psi_{s,n}(M_{n+1})-\psi_{s,n}(M_n)]-\sum_{k=1}^{D-1}\beta_{k,D}(\widetilde{f}_s)\mathbb{E}
\big[M_{n,(k+1);I_{s,n}}^{(D)}-M_{n,(k);I_{s,n}}^{(D)}\big]  \le C_q.
\ee
\end{lem}

\subsection{Proof of Theorem \ref{Theo:Clustering} assuming Lemma \ref{Lem:MainIngrClust}}

Before proving Lemma \ref{Lem:MainIngrClust}, we demonstrate how it implies Theorem \ref{Theo:Clustering}.
Since $\supp(q)$ is bounded we have 
\[
|\varphi_n(v)-\varphi_{n+1}(v)|\le C_q,~~ \hbox{a.s. for $v$ with $|v| = n$.}
\]
This property can be propagated up the tree, 
to see that
\[
|\varphi_n(o)-\varphi_{n+1}(o)|\le C_q,~~\hbox{a.s.}
\]
In other words, we have 
\be\label{aug2902}
|M_{n+1}-M_n|\le C_q, ~~\hbox{a.s.}
\ee

Now, let $s\in \{1,\dots, N_f+1\}$ be arbitrary. The function $\psi_{s,n}$, defined in (\ref{aug2914}), 
is Lipschitz with the Lipschitz constant equal 1. Thus,~(\ref{aug2902}) 
implies 
\begin{equation}\label{eq:PsiBound}
\mathbb{E}|\psi_{s,n}(M_{n+1})-\psi_{s,n}(M_n)| \le C_q . 
\end{equation}
Combining \eqref{eq:PsiBound} with Lemma \ref{Lem:MainIngrClust} yields that for all $D\ge d_0$ we have
\begin{equation}\label{eq:OrdDiffsBound}
-2C_q\le\sum_{k=1}^{D-1}\beta_{k,D}(\widetilde{f}_s)\mathbb{E}\Big[M_{n,(k+1);I_{s,n}}^{(D)}-M_{n,(k);I_{s,n}}^{(D)}\Big] \le 2C_q.  
\end{equation}
By (\ref{aug2904}), we know that  $\widetilde{f}_s$ has no sign change in $[0,1]$. Thus, Proposition \ref{Lem:BernCoeffPos} 
can be applied to either $\widetilde f_s$ or~$(-\widetilde{f}_s)$.
Hence, there is $D_s\ge d_0$ such that the coefficients $(\beta_{k,D_s}(\widetilde{f}_s))$, $k\in\{1,\dots, D_s-1\}$  all have the same  sign  
and at least one of them is not zero.  
Fix some $k_s$ such that $\beta_{k_s,D_s}(\widetilde{f}_s)\neq 0$. Since all $\beta_{k,D_s}(\widetilde{f}_s)$ have the same sign
and because of (\ref{aug2912}),~\eqref{eq:OrdDiffsBound} implies that
\begin{equation}
\mathbb{E}\left[M_{n,(k_s+1);I_{s,n}}^{(D_s)}-M_{n,(k_s);I_{s,n}}^{(D_s)}\right]\le 2C_q|\beta_{k_s,D_s}(\widetilde{f}_s)|^{-1}<\infty \label{eq:BoundIncr}
\end{equation}
is bounded uniformly in $n$. 

It remains to show that \eqref{eq:BoundIncr} implies tightness of $(M_{n;I_s}-\mathrm{med}(M_{n;I_s}))_{n\in\N}$. 
Let us fix $\varepsilon>0$, and denote by~$q_{n,\varepsilon}$ the $\varepsilon$-quantile of $M_{n;I_s}$ so that 
$q_{n,1/2}$ is the median of $M_{n;I_s}$. We have
\be\label{aug2916bis}
\bal
&\mathbb{E}\Big[M_{n,(k_s+1);I_{s,n}}^{(D_s)}-M_{n,(k_s);I_{s,n}}^{(D_s)}\Big] 
\ge (q_{n,1/2}-q_{n,\eps})\Pm\Big[M_{n,(k_s);I_s}\le q_{n,\eps},~M_{n,(k_s+1);I_s}\ge q_{n,1/2}\Big]
\\
&=\binom{D_s}{k_s}\varepsilon^{k_s}(1/2)^{D_s-k_s}(q_{n,1/2}-q_{n,\varepsilon}).
\enbal
\ee
Combining (\ref{aug2916bis}) 
with \eqref{eq:BoundIncr} yields that for all $n\in\N$ we have
\begin{equation}
q_{n,1/2}-q_{n,\varepsilon}\le \left(\binom{D_s}{k_s}\varepsilon^{k_s}(1/2)^{D_s-k_s}\right)^{-1}\cdot \left(2C_q|\beta_{k_s,D_s}(\widetilde{f}_s)|^{-1}\right)<\infty. \label{eq:BdLowTail}
\end{equation}
An analogous argument yields that
\begin{equation}
q_{n,(1-\varepsilon)}-q_{n,1/2} \le \left(\binom{D_s}{k_s}(1-\varepsilon)^{D_s-k_s}(1/2)^{k_s}\right)^{-1}\cdot \left(2C_q|\beta_{k_s,D_s}(\widetilde{f}_s)|^{-1}\right)<\infty  \label{eq:BdUpTail}
\end{equation}
Since $\varepsilon>0$ was arbitrary, together \eqref{eq:BdLowTail} and \eqref{eq:BdUpTail} yield that $(M_{n;I_s}-\mathrm{med}(M_{n;I_s}))_{n\in\N}$ 
is tight.~$\Box$

\subsection{An auxiliary lemma}

It is convenient to introduce the notation 
\be\label{aug3006}
B_{k,d}(x) := \sum_{l=k}^d b_{l,d}(x).
\ee
with this we can write the recursion polynomial \eqref{eq:gThresh} as 
\be
g(x) = \sum_{d=1}^{d_0}\sum_{k=1}^d p_d\zeta_{k,d}B_{k,d}(x)
\ee
and similarly for the nonlinearity $f(u)$. To prove Lemma \ref{Lem:MainIngrClust} we need to understand 
how to expand the polynomial $f((\alpha_2-\alpha_1)x+\alpha_1)$ that appears in (\ref{aug2910}) 
as a weighted sum of $B_{k,d}(x)$. This is done in the next lemma. We recall the notation
\be
\binom{d}{l,j,d-j-l}=\farc{d!}{l!j!(d-j-l)!}.
\ee
\begin{lem}\label{Lem:RescaleNonlin}
Fix a random threshold voting model $(p_d,\zeta_{k,d})_{d\le d_0, k\in\{1,\dots, d\}}$ and $\alpha_1<\alpha_2\in [0,1]$. For the associated nonlinearity $f$ we have
\begin{equation}
\begin{aligned}
f((\alpha_2-\alpha_1)x+\alpha_1) &= \sum_{d=1}^{d_0}\sum_{k=1}^d\sum_{l=0}^{k-1}\sum_{m=k-l}^{d-l}p_d\zeta_{k,d} \binom{d}{l,m,d-m-l}
\alpha_1^l(\alpha_2-\alpha_1)^m(1-\alpha_2)^{d-m-l}B_{k-l,m}(x)\\
&-(\alpha_2-\alpha_1)x+f(\alpha_1).
\end{aligned} \label{eq:RescaleNonlin}
\end{equation}
\end{lem}

\begin{proof}
We recall that, by definition,
\be\label{aug3002}
f((\alpha_2-\alpha_1)x+\alpha_1) = \sum_{d=1}^{d_0}\sum_{k=1}^d p_d\zeta_{k,d} B_{k,d}((\alpha_2-\alpha_1)x+\alpha_1)-(\alpha_2-\alpha_1)x-\alpha_1.
\ee
Thus, to prove \eqref{eq:RescaleNonlin} it is enough to show that for all $d\in\N$, $k\le d$ we have 
\begin{equation}
\begin{aligned}
&B_{k,d}((\alpha_2-\alpha_1)x+\alpha_1)= \sum_{l=0}^{k-1}\sum_{m=k-l}^{d-l} 
\!\binom{d}{l,m,d-m-l}\alpha_1^l(\alpha_2-\alpha_1)^m(1-\alpha_2)^{d-m-l}B_{k-l,m}(x)+B_{k,d}(\alpha_1). 
\end{aligned} \label{eq:fml3}
\end{equation}
We do this by a direct computation. We start by looking at the stretched version of $b_{j,d}$:
\be\label{aug3004}
\bal
&b_{j,d}((\alpha_2-\alpha_1)x+\alpha_1)= \binom{d}{j}((\alpha_2-\alpha_1)x+\alpha_1)^j(1-(\alpha_2-\alpha_1)x-\alpha_1)^{d-j}\\
&=\binom{d}{j} \Big[\sum_{l=0}^j \binom{j}{l} (\alpha_2-\alpha_1)^l x^l\alpha_1^{j-l}\Big]
\Big[\sum_{m=0}^{d-j}\binom{d-j}{m}(1-\alpha_2)^{d-j-m}(\alpha_2-\alpha_1)^{m}(1-x)^{m}\Big]\\
&= \sum_{l=0}^j \sum_{m=0}^{d-j} \binom{d}{j}\binom{j}{l}\binom{d-j}{m} (\alpha_2-\alpha_1)^{l+m}(1-\alpha_2)^{d-j-m}\alpha_1^{j-l}x^l(1-x)^{m}\\
&= \sum_{l=0}^j\sum_{m=0}^{d-j} \binom{d}{j-l,l+m,d-j-m} b_{l,l+m}(x)(\alpha_2-\alpha_1)^{l+m}(1-\alpha_2)^{d-j-m}\alpha_1^{j-l}\\
&=\sum_{l=0}^j\sum_{m=0}^{d-j} \binom{d}{l,j-l+m,d-j-m} b_{j-l,j-l+m}(x)(\alpha_2-\alpha_1)^{j-l+m}(1-\alpha_2)^{d-j-m}\alpha_1^{l}.
\enbal
\ee
The next to last step above uses the definition of $b_{l,l+m}$ as well as the relation 
\be
\bal
\binom{d}{j}\binom{j}{l}\binom{d-j}{m}\Big[\binom{l+m}{l}\Big]^{-1} &=\farc{d!}{j!(d-j)!}\frac{j!}{l!(j-l)!}\farc{(d-j)!}{m!(d-j-m)!}\farc{l!m!}{(l+m)!}\\
&=\farc{d!}{(j-l)!(d-j-m)!(l+m)!}
= \binom{d}{j-l,l+m,d-j-m}.
\enbal
\ee 
Plugging (\ref{aug3004}) into the definition (\ref{aug3006}) of $B_{k,d}$ yields
\be \label{aug3008}
\bal
&B_{k,d}((\alpha_2-\alpha_1)x+\alpha_1) = \sum_{j= k}^d b_{j,d}((\alpha_2-\alpha_1)(x)+\alpha_1)\\
&=\sum_{j=k}^d \sum_{l=0}^j \sum_{m=0}^{d-j} \binom{d}{l,j-l+m,d-j-m}b_{j-l,j-l+m}(x)(\alpha_2-\alpha_1)^{j-l+m}(1-\alpha_2)^{d-j-m}\alpha_1^{l}.
\enbal
\ee
Exchanging the order of summation of the sums over $j$ and $l$ and changing the index of summation of the third sum to $\widehat{m} = j-l+m$ yields
\be
\bal
&B_{k,d}((\alpha_2-\alpha_1)x+\alpha_1)
=\sum_{l=k+1}^d  \sum_{j=l}^d \sum_{m=j-l}^{d-l} \binom{d}{l,m,d-m-l}b_{j-l,m}(x)(\alpha_2-\alpha_1)^{m}(1-\alpha_2)^{d-m-l}\alpha_1^{l}\\
&\qquad+\sum_{l=0}^{k}\sum_{j=k}^d \sum_{m = j-l}^{d-l}\binom{d}{l,m,d-m-l}b_{j-l,m}(x)(\alpha_2-\alpha_1)^{m}(1-\alpha_2)^{d-m-l}\alpha_1^{l}.
\enbal
\ee
Exchanging the order of summation of the sums over $m$ and  $j$ yields
\be 
\bal
B_{k,d}((\alpha_2-\alpha_1)x+\alpha_1)&=
\sum_{l=k+1}^d   \sum_{m=0}^{d-l} \sum_{j=l}^{m+l} \binom{d}{l,m,d-m-l}b_{j-l,m}(x)(\alpha_2-\alpha_1)^{m}(1-\alpha_2)^{d-m-l}\alpha_1^{l}\\
&+\sum_{l=0}^{k} \sum_{m = k-l}^{d-l} \sum_{j=k}^{m+l}\binom{d}{l,m,d-m-l}b_{j-l,m}(x)(\alpha_2-\alpha_1)^{m}(1-\alpha_2)^{d-m-l}\alpha_1^{l}.
\enbal
\ee
Now, switching the summation over $j$ to $\hat j=j-l$ and dropping the hat gives 
\be\label{aug3016}
\bal
B_{k,d}((\alpha_2-\alpha_1)x+\alpha_1)&=\sum_{l=k+1}^d   \sum_{m=0}^{d-l} \sum_{j=0}^{m} \binom{d}{l,m,d-m-l}b_{j,m}(x)(\alpha_2-\alpha_1)^{m}(1-\alpha_2)^{d-m-l}\alpha_1^{l}\\
&+\sum_{l=0}^{k} \sum_{m = k-l}^{d-l} \sum_{j=k-l}^{m}\binom{d}{l,m,d-m-l}b_{j,m}(x)(\alpha_2-\alpha_1)^{m}(1-\alpha_2)^{d-m-l}\alpha_1^{l}.
\enbal
\ee
Next, we use \eqref{Def:BernSum} in the first sum in (\ref{aug3016}) and the definition (\ref{aug3006}) of $B_{k-l,m}$ in the second, to obtain
\be \label{aug3018}
\bal
B_{k,d}((\alpha_2-\alpha_1)x+\alpha_1)&= \sum_{l=k+1}^d \sum_{m=0}^{d-l} \binom{d}{l,m,d-m-l}(\alpha_2-\alpha_1)^m(1-\alpha_2)^{d-m-l}\alpha_1^l\\
&+ \sum_{l=0}^{k}\sum_{m=k-l}^{d-l} \binom{d}{l,m,d-m-l} \alpha_1^l(\alpha_2-\alpha_1)^m(1-\alpha_2)^{d-m-l} B_{k-l,m}(x).
\enbal
\ee
The summation over $m$ in the  first sum in (\ref{aug3018}) can be re-written as
\be\label{aug3020}
\bal
&\sum_{m=0}^{d-l} \binom{d}{l,m,d-m-l}(\alpha_2-\alpha_1)^m(1-\alpha_2)^{d-m-l}=
\sum_{m=0}^{d-l}\farc{d!(d-l)!}{l!m!(d-m-l)!(d-l)!}(\alpha_2-\alpha_1)^m(1-\alpha_2)^{d-m-l}\\
&=\farc{d!}{l!(d-l)!}(1-\alpha_1)^{d-l}=\binom{d}{l}(1-\alpha_1)^{d-l}.
\enbal
\ee
Furthermore, as $B_{0,m}(x)\equiv 1$ because of (\ref{Def:BernSum}), the $l=k$ summand in the second sum in (\ref{aug3018}) can be written as 
\be\label{aug3021}
\bal
&\sum_{m=0}^{d-k} \binom{d}{k,m,d-m-k} \alpha_1^k(\alpha_2-\alpha_1)^m(1-\alpha_2)^{d-m-k} B_{0,m}(x)\\
&=
\sum_{m=0}^{d-k}\farc{d!(d-k)!}{k!m!(d-m-k)!(d-k)!} \alpha_1^k(\alpha_2-\alpha_1)^m(1-\alpha_2)^{d-m-k} =\binom{d}{k}\alpha_1^k(1-\alpha_1)^{d-k}.
\enbal
\ee
Using (\ref{aug3020}) and (\ref{aug3021}) in (\ref{aug3018}) leads to 
\begin{align*}
&B_{k,d}((\alpha_2-\alpha_1)x+\alpha_1)\\
&= \sum_{l=k}^d \binom{d}{l} \alpha_1^l (1-\alpha_1)^{d-l}+\sum_{l=0}^{k-1}\sum_{j=k-l}^{d-l} \binom{d}{l,j,d-j-l} \alpha_1^l(\alpha_2-\alpha_1)^j(1-\alpha_2)^{d-j-l}B_{k-l,j}(x)\\
&= B_{k,d}(\alpha_1)+\sum_{l=0}^{k-1}\sum_{m=k-l}^{d-l} \binom{d}{l,m,d-m-l} \alpha_1^l(\alpha_2-\alpha_1)^m(1-\alpha_2)^{d-m-l}B_{k-l,m}(x).
\end{align*}
This proves \eqref{eq:fml3} finishing the proof of \eqref{eq:RescaleNonlin}.
\end{proof}

\subsection{Proof of Lemma \ref{Lem:MainIngrClust}}

Here, we prove Lemma \ref{Lem:MainIngrClust}, finishing the proof to Theorem \ref{Theo:Clustering}. Let us consider 
a collection of independent random variables  $(Z, (M_{n,k})_{k\in\N}$, $(X_k)_{k\in\N})$ such that
\be\label{aug3024}
\bal
&\mathbb{P}[Z= d] = p_d,\\
&M_{n,k}\stackrel{d}{=}M_n,\\
&X_k\sim q,
\enbal
\ee
and also a random variable $L$ independent of $(M_{n,k}, X_k)_{k\in\N}$ with 
\be
\mathbb{P}[L = k| Z=d] = \zeta_{k,d}.
\ee 
Observe that the same reasoning
giving the recursion \eqref{eq:RecEqThresh} yields a recursion relation 
\begin{equation}\label{eq:RecursionMn}
M_{n+1} \stackrel{d}{=}\sum_{d=1}^{d_0} \one({Z=d})(M_{n,1}+X_1,\dots, M_{n,d}+X_d)_{(L)}.  
\end{equation}

Recall that the Lipschitz function $\psi_{s,n}(x)$ that appears in the statement of Lemma~\ref{Lem:MainIngrClust} is defined by~(\ref{aug2914}).
Since there is only one non-zero term in the sum in the right side of (\ref{eq:RecursionMn}), 
the recursion \eqref{eq:RecursionMn} immediately implies that
\be\label{aug3022}
\bal
\mathbb{E}[\psi_{s,n}(M_{n+1})] &=\mathbb{E}\Big[\psi_{s,n}\Big(\sum_{d=1}^{d_0}\one({Z = d})(M_{n,1}+X_1,\dots, M_{n.d}+X_d)_{(L)}\Big)\Big]\\
&= \mathbb{E}\Big[\sum_{d=1}^{d_0}\one({Z=d})\psi_{s,n}\left((M_{n,1}+X_1,\dots, M_{n,d}+X_d)_{(L)}\right)\Big]\\
&\le \mathbb{E}\Big[\sum_{d=1}^{d_0}\one({Z=d})\psi_{s,n}\left(M_{n,(L)}^{(d)}\right)\Big]\!+\mathbb{E}\Big[\max_{k\in \{1,\dots, d_0\}} X_k\Big]
\le \sum_{d=1}^{d_0}\sum_{k=1}^d p_d\zeta_{k,d}\mathbb{E}[\psi_{s,n}(M_{n,(k)}^{(d)})]+C_q.
\enbal
\ee
The third step used the fact that $\psi_{s,n}$ is Lipschitz with Lipschitz constant $1$ and the last step used (\ref{aug3025}) and (\ref{aug3024}). 

Let us write 
\be\label{aug3028}
\psi_{s,n} (M_n)= (q_{s,n}-q_{s-1,n})\one(M_n\le q_{s-1,n}) +(q_{s,n}-M_n)\one({M_n\in [q_{s-1,n},q_{s,n}]}) .
\ee
and take the expectation:
\be\label{aug3027}
\mathbb{E}[\psi_{s,n}(M_n)] = \alpha_{s-1}(q_{s,n}-q_{s-1,n})+(\alpha_s-\alpha_{s-1})\mathbb{E}[q_{s,n}-M_{n;I_{s,n}}].
\ee
An analogous argument using $\min_{k\in\{1,\dots, d_0\}} X_k$ and subtracting (\ref{aug3027})  
from both sides of (\ref{aug3022}) yields 
\begin{equation}\label{eq:DefSigmasn}
\begin{aligned}
&\mathbb{E}[\psi_{s,n}(M_{n+1})-\psi_{s,n}(M_n)]-C_q\\
&\le \sum_{d=1}^{d_0}\sum_{k=1}^d p_d\zeta_{k,d}\mathbb{E}[\psi_{s,n}(M_{n,(k)}^{(d)})]-\alpha_{s-1}(q_{s,n}-q_{s-1,n})
-(\alpha_s-\alpha_{s-1})\mathbb{E}[q_{s,n}-M_{n;I_s,n}]  \\
&\le  \mathbb{E}[\psi_{s,n}(M_{n+1})-\psi_{s,n}(M_n)]+C_q.
\end{aligned}  
\end{equation}
By decomposing with regard to how many of the $M_{n,k}$ are in $(-\infty, q_{s-1,n}]$, in $[q_{s-1,n}, q_{s,n}]$ and in $[q_{s,n},\infty)$,
respectively, we get from (\ref{aug3028})
\be\label{aug3029}
\bal
\mathbb{E}\left[\psi_{s,n}\left(M_{n,(k)}^{(d)}\right) \right] &= \sum_{l=0}^{k-1}\sum_{j=k-l}^{d-l} \binom{d}{l,j,d-j-l}\alpha_{s-1}^l(\alpha_s-\alpha_{s-1})^{j}(1-\alpha_{s})^{d-l-j}\mathbb{E}[q_{s,n}-M_{n,(k-l); I_{s,n}}^{(j)} ]\\
& +\sum_{l=k}^d \sum_{j=0}^{d-l} \binom{d}{l,j,d-j-l}\alpha_{s-1}^l(\alpha_{s}-\alpha_{s-1})^{j}(1-\alpha_{s})^{d-l-j}(q_{s,n}-q_{s-1,n} ).
\enbal
\ee
The second sum in the right side can be simplified by writing
\be\label{aug3030}
\bal
&\sum_{j=0}^{d-l} \binom{d}{l,j,d-j-l}(\alpha_{s}-\alpha_{s-1})^{j}(1-\alpha_{s})^{d-l-j}\\
&=\sum_{j=0}^{d-l}\farc{d!(d-l)!}{l!j!(d-j-l)!(d-l)!}(\alpha_{s}-\alpha_{s-1})^{j}(1-\alpha_{s})^{d-l-j}=\binom{d}{l}(1-\alpha_{s-1})^{d-l}.
\enbal
\ee
Let us consider the term in the second line of (\ref{eq:DefSigmasn}):
\be\label{aug3033}
\Sigma_{{s,n}}:=\sum_{d=1}^{d_0}\sum_{k=1}^d p_d\zeta_{k,d}\mathbb{E}[\psi_{s,n}(M_{n,(k)}^{(d)})]-\alpha_{s-1}(q_{s,n}-q_{s-1,n})
-(\alpha_s-\alpha_{s-1})\mathbb{E}[q_{s,n}-M_{n;I_s,n}] .  
\ee
Note that  (\ref{eq:DefSigmasn}) says  
\be\label{aug3032}
\Big|\mathbb{E}[\psi_{s,n}(M_{n+1})-\psi_{s,n}(M_n)]-\Sigma_{s,n}\Big|\le C_q.
\ee
Thus,  the conclusion of Lemma~\ref{Lem:MainIngrClust} will follow if we show that 
\be
\Sigma_{s,n}=\sum_{k=1}^{D-1}\beta_{k,D}(\widetilde{f}_s)\mathbb{E}
\big[M_{n,(k+1);I_{s,n}}^{(D)}-M_{n,(k);I_{s,n}}^{(D)}\big].
\ee
Using (\ref{aug3029}) and (\ref{aug3030}) in the definition (\ref{aug3033}) of $\Sigma_{s,n}$, we can re-write that sum as 
\begin{equation}\label{eq:Sigmasnsimpl}
\begin{aligned}
\Sigma_{s,n} &=\sum_{d=1}^{d_0}\sum_{k=1}^d\sum_{l=0}^{k-1}\sum_{j=k-l}^{d-l}p_d\zeta_{k,d}\binom{d}{l,j,d-j-l}\alpha_{s-1}^l(\alpha_s-\alpha_{s-1})^j(1-\alpha_s)^{d-l-j}\mathbb{E}[q_{s,n}-M_{n,(k-l);I_{s,n}}^{(j)} ]\\
& +(q_{s,n}-q_{s-1,n})\Big(\sum_{d=1}^{d_0}\sum_{k=1}^dp_d\zeta_{k,d}\sum_{l=k}^d\binom{d}{l}\alpha_{s-1}^l(1-\alpha_{s-1})^{d-l}-\alpha_{s-1}\Big)\\
& -(\alpha_s-\alpha_{s-1})\mathbb{E}[q_{s,n}-M_{n;I_{s,n}}]\\
&=\sum_{d=1}^{d_0}\sum_{k=1}^d\sum_{l=0}^{k-1}\sum_{j=k-l}^{d-l}p_d\zeta_{k,d} \binom{d}{l,j,d-j-l}\alpha_{s-1}^l(\alpha_s-\alpha_{s-1})^j(1-\alpha_s)^{d-l-j}\mathbb{E}[q_{s,n}-M_{n,(k-l);I_{s,n}}^{(j)}]\\
&\qquad-(\alpha_s-\alpha_{s-1})\mathbb{E}[q_{s,n}-M_{n;I_{s,n}}].
\end{aligned}  
\end{equation}
The second step above used the identity
\[
\sum_{d=1}^{d_0}\sum_{k=1}^dp_d\zeta_{k,d}\sum_{l=k}^d\binom{d}{l}\alpha_{s-1}^l(1-\alpha_{s-1})^{d-l}-\alpha_{s-1}  = f(\alpha_{s-1}) = 0.
\]
On the other hand, by using Lemma \ref{Lem:RescaleNonlin} as well as $f(\alpha_{s-1}) = 0$ we see that 
\begin{equation}
\begin{aligned}
&\widetilde{f}_s(x)= \sum_{d=1}^{d_0}\sum_{k=1}^d\sum_{l=0}^{k-1}\sum_{j=k-l}^{d-l}p_d\zeta_{k,d}\binom{d}{l,j,d-j-l}\alpha_{s-1}^l(\alpha_s-\alpha_{s-1})^j(1-\alpha_s)^{d-j-l}B_{k-l,j}(x)\\
&\qquad-(\alpha_s-\alpha_{s-1})x. 
\end{aligned}\label{eq:fstildExp}
\end{equation}
Comparing \eqref{eq:Sigmasnsimpl} to \eqref{eq:fstildExp}  we see that the coefficient in front of $B_{k-l,j}(x)$ in the expression for
$\widetilde{f}_s$ equals the coefficient in front of $\mathbb{E}[q_{s,n}-M_{n,(k-l);I_{s,n}}^{(j)}]$ in $\Sigma_{s,n}$. Thus to show that for all $D\ge d_0$ there are coefficients $\beta_{k,D,s}$ such that
\begin{equation}
\Sigma_{s,n} =  \sum_{k=1}^D \beta_{k,D,s}\mathbb{E}[q_{s,n}-M_{n,(k);I_{s,n}}^{(D)}]\quad\text{and}\quad \widetilde{f}_s(x) = \sum_{k=1}^D\beta_{k,D,s}B_{k,D}(x) \label{eq:CoeffSigftildequal}
\end{equation}
it is enough to show that there is a family of multi-linear functions 
\[
f_{k,d,D}:\R^D\to \R,~~ D\ge 1, ~~d\in \{1,\dots, D\}, ~~k\in\{1,\dots, d\},
\]
such that
\begin{equation}
\begin{aligned}
B_{k,d} &= f_{k,d,D}(B_{1,D},\dots, B_{D,D}),
\enbal
\ee
and
\be
\bal
\mathbb{E}[q_{s,n}-M_{n,(k);I_{s,n}}^{(d)}] &= f_{k,d,D}\left(\mathbb{E}[q_{s,n}-M_{n,(1);I_{s,n}}^{(D)}],\dots, \mathbb{E}[q_{s,n}-M_{n,(D);I_{s,n}}^{(D)}]\right).
\end{aligned}
\end{equation}
Since a composition of multilinear functions is multilinear itself, it is enough to show this for $D = d+1$. 

To this end, first note that, using \eqref{eq:BernDegEl} and the definition of $B_{k,d}$, we have for all $d\in\N$, $k\le d$
\be\label{eq:Bkd1Step}
\bal
&B_{k,d}(x) = \sum_{j=k}^d b_{j,d}(x) = \sum_{j=k}^d \Big[\frac{d-j+1}{d+1}b_{j,d+1}(x)+\frac{j+1}{d+1}b_{j+1,d+1}(x)\Big]\\
&=\sum_{j=k}^d \Big[b_{j,d+1}(x)-\frac{j}{d+1}b_{j,d+1}(x)+\frac{j+1}{d+1}b_{j+1,d+1}(x)\Big]
=\sum_{j=k}^d b_{j,d+1}(x)-\frac{k}{d+1}b_{k,d+1}(x)+b_{d+1,d+1}(x)\\
&=B_{k,d+1}(x)-\frac{k}{d+1}[B_{k,d+1}(x)-B_{k+1,d+1}(x)] 
=\frac{d+1-k}{d+1}B_{k,d+1}(x)+\frac{k}{d+1}B_{k+1,d+1}(x).
\enbal
\ee
On the other hand, for any collection of \iid random variables $X_k$, $k\in\Nm$, so that $X_1$ has a continuous density, and any $d\ge 1$, $1\le k\le d$, we have the identity
\be\label{aug3104}
\Em[X^{(d)}_{(k)}]=\frac{k}{d+1}\Em[X^{(d+1)}_{(k+1)}]+\frac{d+1-k}{d+1}\Em[X^{(d+1)}_{(k)}],
\ee 
as can be seen simply by adding $X_{d+1}$ to the collection $\{X_1,\dots,X_d\}$ and looking at whether $X_{d+1}$ is to the left or to the right
of $X_{(k)}^{(d)}$.
Applying (\ref{aug3104}) to  $M_{n,k;I_{s,n}}$  shows that for all $d\in\N$, $k\le d$ we have
\begin{equation}
\mathbb{E}[q_{s,n}-M_{n,(k);I_{s,n}}^{(d)}] = \frac{k}{d+1}\mathbb{E}[q_{s,n}-M_{n,(k+1);I_{s,n}}^{(d+1)}]+\frac{d+1-k}{d+1}\mathbb{E}[q_{s,n}-M_{n,(k);I_{s,n}}^{(d+1)}].\label{eq:Mnk1Step}
\end{equation}

As mentioned above, comparing \eqref{eq:Bkd1Step} to \eqref{eq:Mnk1Step} yields that for all $D\ge d_0$ there are coefficients $\beta_{k,D,s}$ such that \eqref{eq:CoeffSigftildequal} holds. As we also have
 \[
 b_{k,D} = B_{k,D}-B_{k+1,D},
 \]
and 
\[
\mathbb{E}[M_{n,(k+1);I_{s,n}}^{(D)}-M_{n,(k);I_{s,n}}^{(D)}] = \mathbb{E}[q_{s,n}-M_{n,(k);I_{s,n}}^{(D)}]-\mathbb{E}[q_{s,n}-M_{n,(k+1);I_{s,n}}^{(D)}],
\]
equation \eqref{eq:CoeffSigftildequal} implies that for all $s\in \{1,\dots, N_f+1\}$ and $D\ge d_0$ we have
\begin{equation}
\Sigma_{s,n} = \sum_{k=1}^{D-1}\beta_{k,D}(\widetilde{f}_s)\mathbb{E}\left[M_{n,(k+1);I_{s,n}}^{(D)}-M_{n,(k);I_{s,n}}^{(D)}\right]. \label{eq:SigmasnLastStep}
\end{equation}
Recalling  \eqref{eq:DefSigmasn} we see that \eqref{eq:SigmasnLastStep} implies the statement of Lemma \ref{Lem:MainIngrClust}, finishing the proof.
$\Box$

\subsection{The binary-ternary case as an example}
\label{subsec-bt}

In this section we will look at the threshold voting models with $p_2 = p$, $p_3 = (1-p)$ for $p\in [0,1]$ and~$\zeta_{2,3} = \zeta_{2,2} = 1$. 
In other words, a parent who has three children is assigned the middle one of their values, while a parent with two children gets the larger value of the two.

There are several 
reasons to look at these models: they have an additional probabilistic interpretation, they are convenient for showing that we can get a slightly stronger result than 
Theorem~\ref{Theo:Clustering} with probabilistic means, and for $p<1/2$ they are example for which the nonlinearity has the single additional zero $1/(2(1-p))$ in $(0,1)$ as 
well as $f'(0)<0$, $f'(1)<0$. In Section \ref{Sec:TightnessAnMeans}, we will use analytic methods to show that for such nonlinearities $f$ 
the sequence $(M_n-\mathrm{med}(M_n))_{n\in\N}$ is tight.

Let us mention an alternative probabilistic interpretation for that voting model.
Let~$\cT_n$ be the genealogical tree of the underlying BRW up to generation $n$ 
and~$\mathbf{T}_n^{(2)} := \{(T_n,o)\subseteq \cT_n : T_n\ \text{binary}\}$   be the collection of rooted full binary subtrees of $\cT_n$ with root~$o$ and depth $n$. Given a binary subtree~$(T_n,o)\in \mathbf{T}_n^{(2)}$,  we define
\[
M_{n,T_n} := \max_{v\in T_n : |v| = n} S_v,
\]
as the maximum at time $n$ along $T_n$. Finally, we set
\begin{equation}
  \label{eq-1711d}
\widetilde{M}_n := \min_{(T_n,o)\in \mathbf{T}_n^{(2)}} M_{n,T_n}
\end{equation}
to be the smallest maximum along all binary subtrees of $\cT_n$. Analogously, we define 
\begin{equation}
  \label{eq-1711e}
\widehat{M}_n := \max_{(T_n,o)\in \mathbf{T}_n^{(2)}} \min_{v\in T_n : |v| =  n} S_v
\end{equation} 
as the largest minimum along all binary subtrees.

Let us make a couple of simple observations.
First, it follows from the definition above that
\be\label{aug3110}
M_n = \widetilde{M}_n.
\ee
Furthermore, in the case of purely ternary branching $p=0$, we have
\be\label{aug3112}
M_n = \widetilde{M}_n = \widehat{M}_n.
\ee
In particular, it follows that, for $q$ symmetric and purely ternary branching, 
the distribution of $M_n$ is symmetric for all $n\ge 1$, and 
\be\label{aug3114}
\mathbb{E}[M_n] = \mathrm{med}(M_n) = 0.
\ee

While the description of $M_n$ as the smallest maximum of all binary subtrees of $\cT_n$ is quite nice and links the study of $M_n$ to the study of the maximum of BRWs,   
we were unable to use it to gain any insights into the distribution of $M_n$. One of the reasons for
this is that while we have very precise control of $\mathbb{P}[M_{n,T_n}\le t]
$ for $T_n$ a fixed binary subtree of $\cT_n$, there are $3^{2^n-1}$ binary subtrees of $\cT_n$, which are far too many for a first moment method to work. Of course, many of 
these binary subtrees share many vertices. For example, for any given binary subtree $T_n$ of $\cT_n$ there are at least $3^{2^{n-1}}$ binary subtrees $\widetilde{T}_n$ such 
that~$\{v\in T_n : |v| = n-1\} = \{v\in \widetilde{T}_n : |v| = n-1\}$. The issue we could not overcome is that we do not know how to properly use the fact that many of the 
maxima along the binary subtrees are strongly correlated.

The rest of this section is devoted to the 
subsequential tightness of $(M_n)_{n\in\N}$ in the fully ternary case.
\begin{thm}\label{Theo:SubseqTight}
In addition to our standing assumptions (\ref{aug3025}) on $q$, assume that $q$ is symmetric. 
Consider the voting model with $p_3 = 1$ and $\zeta_{2,3} = 1$. There is a subset $I\subseteq \N$ of the natural numbers such 
that~$(M_n)_{n\in I}$ is tight and 
\be\label{sep102}
\liminf_{n\to \infty}\frac{I\cap \{1\dots, n\}}{n}>0.
\ee
\end{thm}
Let us first outline the proof of this theorem.  It relies on the voting
model interpretation of $M_n$. 
Observe that in the symmetric case we have 
\[
-M_{n;I_1} \stackrel{d}{=} M_{n;I_2} \stackrel{d}{=} |M_n|.
\]
Thus, Theorem \ref{Theo:Clustering} implies that $(|M_n|-\mathrm{med}(|M_n|))_{n\in\N}$ is tight. Thus, for the full tightness of $M_n$ it is enough to show that there 
are $\varepsilon, C>0$ such that for all $n\in\N$ we have 
\be\label{sep104}
\mathbb{P}[|M_n|\le C]\ge \varepsilon.
\ee
We have not been able to prove (\ref{sep104}) by purely probabilistic means. 
Instead we show by contradiction that if 
\be\label{sep106}
\liminf_{n\to\infty}\big[\mathrm{med}(|M_n|)\big]=+\infty,
\ee
it is too likely for a particle $v$ with $S_v\approx 0$ to be voted to the top, 
making $\mathbb{P}[M_n\approx 0]$ too large. 
To see this, first note that $S_v$  is voted to the top if and only if
at each ancestor $v_k$, $|v_k| = k$, of $v$ we have
\begin{equation}
\max_{w\in D_1(v_k)\setminus \{v_{k+1}\}} \varphi_n(w)\ge S_v\ge \min_{w\in D_1(v_k)\setminus \{v_{k+1}\}} \varphi_n(w). \label{eq:KeySubseqTight}
\end{equation}
Suppose now that (\ref{sep106}) holds and take $N$ sufficiently large. Because of (\ref{sep106}), we have
\be\label{sep114}
\mathrm{med}(|M_n|)\gg 100N,
\ee
for all $n$ sufficiently large. 
There exists $\eta_N>0$ so that 
with the probability  $(1-\eta_N)^n$, we have, for all~$k\le n$, both
 \be\label{sep110}
 \max_{w\in D_1(v_k)\setminus \{v_{k+1}\}} |S_w-S_{v_k}|\le N,
 \ee
 and 
\be\label{sep108}
\max_{k\le n} |S_v-S_{v_k}|\le N.
\ee
We may also choose $\eta_N$ so that 
\be
\eta_N\to 0,~~\hbox{as $N\to+\infty$.}
\ee
Note that, under the conditions (\ref{sep110}) and (\ref{sep108}),  \eqref{eq:KeySubseqTight} holds if 
 \begin{equation}
 \max_{w\in D_1(v_k)\setminus\{v_{k+1\}}} [\varphi_n(w)-S_w]\ge 2N\quad\text{and} \min_{w\in D_1(v_k)\setminus\{v_{k+1\}}} [\varphi_n(w)-S_w]\le -2N. \label{eq:Comparis}
 \end{equation}
By construction, we have 
\be\label{sep116}
\varphi_n(w)-S_w\stackrel{d}{=} M_{n-|w|}.
\ee
If (\ref{sep114}) holds for all $k\le n$,
the tightness of $(|M_n|-\mathrm{med}(|M_n|))_{n\in\N}$ and (\ref{sep116})
ensure that the probability of the event in \eqref{eq:Comparis} is roughly equal to $2(1/2-\varepsilon_N)^2$ with 
\be
\lim_{N\to+\infty}\eps_N=0.
\ee
Thus, overall we have
\[
\mathbb{P}[M_n\in [-N,N]] = 3^n\mathbb{P}[S_v\in [-N,N], M_n = S_v] \approx 3^n(1-\eta_N)^n2^n (1/2-\varepsilon_N)^{2n},
\]
which is bigger than $1$ for $\eta_N$, $\varepsilon_N$ small enough, which yields a contradiction to (\ref{sep114}). In the actual proof of Theorem \ref{Theo:SubseqTight} 
we will need to strengthen the lower bound so that it still holds (and is bigger than 1) if every once in a while we do not have $\mathrm{med}(|M_{k}|)\gg 100N$.

\subsubsection*{Proof of Theorem \ref{Theo:SubseqTight}}

Given 
\be\label{sep302}
K>C_q,
\ee
sufficiently large, we set
\[
I_{K} := \{n\in\N_0 : \mathrm{med}(|M_n|) \le K\}.
\]
We will put further restrictions on $K$, in addition to (\ref{sep302}), during the proof, keeping it as large as needed, but  independent of~$n$. 

As we have  mentioned, the symmetry of $M_n$ and Theorem \ref{Theo:Clustering}  imply that $(|M_n|-\mathrm{med}(|M_n|))_{n\in\N}$ is tight. 
Thus, the family $(M_n)_{n\in I_{4K}}$ is tight and it is enough to prove that  
\be\label{sep304}
\liminf_{n\to\infty} \frac{I_{4K}\cap\{0,\dots, n-1\}}{n}>0.
\ee
First, we note that for $K$ big enough we have 
\be\label{sep306}
\mathbb{P}[M_k\le -3K] = \mathbb{P}[M_k\ge 3K] \ge 1/2-1/200,~~\hbox{ for all $k\in I_{4K}^c$}.
\ee
This is true, since, using the symmetry of $M_n$,  we can write
\be\label{sep308}
\mathbb{P}[M_k\ge 3K] = \frac{1}{2}\mathbb{P}[|M_k|\ge 3K] 
\ge \farc{1}{2}\mathbb{P}\left[\left| |M_k|-\mathrm{med}(|M_k|)\right|\le K\right]\ge \frac{1}{2}-\frac{1}{200}, ~~\hbox{ for all $k\in I_{4K}^c$},
\ee
as long as $K$ is chosen to be large enough, but independent of $n$.
The last step in (\ref{sep308}) used the tightness of $(|M_n|-\mathrm{med}(|M_n|))_{n\in\N}$. 

Next, we set 
\be\label{sep310}
\bal
&\widetilde{I}_{4K,n} := \{k\in\{0,\dots, n-1\} : n-k-1\in I_{4K}\}.
\enbal
\ee
Finally, we fix a vertex $v$ with $|v| = n$ and define the event
\be\label{sep518}
A_{n,\delta_0}(v) := \left\{|S_{v_k}-S_v|\le K \hbox{ for all $k\le n$} \hbox{ and }   |S_v-S_{v_k}|\le \delta_0
\hbox{ for all } k\in \widetilde{I}_{4K,n}\right\}.
\ee
Here, $\delta_0>0$ is chosen so that 
\be\label{sep506}
\int_{(-\delta_0,\delta_0)^c} q(x)\dx>0, 
\inf_{x\in [-\delta_0, \delta_0]} \int_{-\delta_0}^{\delta_0} q(y-x)\dy>0.
\ee
To see that such choice of $\delta_0>0$ is possible, we use the continuity of $q(x)$ to find $\delta_0>0$ such that
\be\label{sep502}
\int_{[-\delta_0,\delta_0]^c}q(y)\dy=\farc{1}{2}.
\ee
As $q(x)$ is symmetric, this is equivalent to 
\be\label{sep504}
\int_{\delta_0}^\infty q(y)\dy=\farc{1}{4},~~\int_0^{\delta_0}q(y)\dy=\farc{1}{4}.
\ee
Note that then for each $x\in[0,\delta_0]$ we have
\be\label{sep508}
\int_{-\delta_0}^{\delta_0}q(y-x)\dy=\int_{-\delta_0-x}^{\delta_0-x}q(y)\dy\ge\int_{-\delta_0}^0q(y)\dy=\farc{1}{4}.
\ee
By symmetry, we also have, for each $x\in[-\delta_0,0]$:
\be\label{sep510}
\int_{-\delta_0}^{\delta_0}q(y-x)\dy=\int_{-\delta_0-x}^{\delta_0-x}q(y)\dy\ge\int_0^{\delta_0}q(y)\dy=\farc{1}{4}.
\ee
Summarizing (\ref{sep502}) and (\ref{sep508})--(\ref{sep510}), we have chosen $\delta_0>0$ such that (\ref{sep506}) holds.

Using the exchangeability of the vertices in the same generation of $\cT_n$ yields that
\begin{align}\label{eq:Step1SubseqTight}
\mathbb{P}[M_n\in [-K,K]] &= \sum_{w\in \cT_n : |w|=n} \mathbb{P}[M_n\in [-K,K], S_w = M_n]= 3^n\mathbb{P}[M_n = S_v, S_v\in [-K,K]] \notag\\
&\ge 3^n\mathbb{P}\left[ A_{n,\delta_0}(v),\ \forall_{k\in \{0,\dots, n-1\}} \varphi_n(v_k) = \varphi_n(v_{k+1})\in [-K,K]\right].  
\end{align} 
For $k\in \{0,\dots, n-1\}$ let $D_1(v_k) = \{v_{k+1}, w_1(k), w_2(k)\}$ denote the direct descendants of $v_k$. 
Using the exchangeability of $w_1(k)$, $w_2(k)$ yields that we can continue from \eqref{eq:Step1SubseqTight} to get
\begin{align}\label{eq:Step2SubseqTight}
\mathbb{P}[M_n\in [-K,K]] &\ge 6^n \mathbb{P}\left[A_{n,\delta_0}(v), \forall_{k\le n-1} \varphi_n(w_1(k))\le \varphi_n(v_k)\le \varphi_n(w_2(k))\right] \notag\\
&= 6^n \mathbb{P}\left[A_{n,\delta_0}(v), \forall_{k\le n-1} \varphi_n(w_1(k))\le S_v\le \varphi_n(w_2(k))\right].
\end{align} 
The last step used that on the event under consideration we have $\varphi_n(v_k) = S_v$ for all $k\le n$.

To bound the right side of \eqref{eq:Step2SubseqTight} from below,
we need to look at $k\in \widetilde{I}_{4K,n}$ and $k\in \widetilde{I}_{4K,n}^c$ separately. 
For this, we consider the increments
\[
X_{i,k} := S_{w_{i}(k)}-S_{v_k} \sim q,~~i\in\{1,2\}.
\] 
First, for $k\in \widetilde{I}_{4K,n}$ we use that on $A_{n,\delta_0}(v)$ we have $|S_v-S_{v_k}|\le \delta_0$ and thus
\begin{equation}\label{eq:InI4C1}
\begin{aligned}
A_{n,\delta_0}(v)\cap \{\varphi_n(w_1(k))\le S_v\}&\supseteq \{\varphi_n(w_1(k))-S_{v_k}\le -\delta_0\}\cap A_{n,\delta_0}(v)\\
&\supseteq A_{n,\delta_0}(v)\cap \{\varphi_n(w_1(k))-S_{w_1(k)}\le 0\}\cap \{X_{1,k}\le -\delta_0\}.
\end{aligned} 
\end{equation}
By symmetry, we also have
\begin{equation} \label{eq:InI4C2}
A_{n,\delta_0}(v)\cap\{\varphi_n(w_2(k))\ge S_v\}\supseteq A_{n,\delta_0}(v)\cap \{\varphi_n(w_2(k))-S_{w_2(k)}\ge 0\}\cap \{X_{2,k}\ge \delta_0\}. 
\end{equation}
Next, for $k\in \widetilde{I}_{4K,n}^c$ we use that on $A_{n,\delta_0}(v)$ we have $S_v\in [-K,K]$ to see that
\be\label{sep316}
A_{n,\delta_0}(v)\cap \{\varphi_n(w_1(k))\le S_v\}\supseteq A_{n,\delta_0}(v)\cap \{\varphi_n(w_1(k))\le -K\}.
\ee
Furthermore, on $A_{n,\delta_0}(v)$ we have, because of (\ref{sep302}): 
\be
|S_{w_i(k)}| = |S_{v_k}+X_{i,k}|\le |S_{v_k}|+|X_{i,k}|\le K+C_q\le 2K.
\ee
This, together with (\ref{sep316}) and, once again (\ref{sep302}), implies  
\begin{equation}
A_{n,\delta_0}(v)\cap \{\varphi_n(w_1(k))\le S_v\}\supseteq A_{n,\delta_0}(v) \cap\{\varphi(w_1(k))-S_{w_1(k)}\le -3K\}. \label{eq:InI4Cc1}
\end{equation}
Another use of symmetry yields the analog of (\ref{eq:InI4Cc1})
\begin{equation}
A_{n,\delta_0}(v)\cap\{\varphi_n(w_2(k))\ge S_v\}\supseteq A_{n,\delta_0}(v)\cap\{\varphi_n(w_2(k))-S_{w_2(k)} \ge 3K\}.\label{eq:InI4Cc2}
\end{equation}

We note that 
\[
\varphi_n(w_i(k))-S_{w_i(k)}\stackrel{d}{=} M_{n-k-1}
\]
depends only on the increments of the descendants of $w_i(k)$, while  $A_{n,\delta_0}(v)$  is measurable with
respect to the increments on the path on the tree~$\cT_n$ 
that connects the vertex~$v$ to the root $o$. 
Thus, the random variables $\one_{A_{n,\delta_0}(v)}$ and~ $(\varphi_n(w_i(k))-S_{w_i(k)})_{k\le n-1, i\in\{1,2\}})$ are independent from each other. 
Using this consideration, together with \eqref{eq:InI4C1} and~\eqref{eq:InI4C2} for   $k\in \widetilde{I}_{4K,n}$, 
as well as~\eqref{eq:InI4Cc1} and~\eqref{eq:InI4Cc2} for $k\in \widetilde{I}_{4K,n}^c$, 
in~\eqref{eq:Step2SubseqTight} yields  
\be 
\bal
\mathbb{P}\big[M_n\in [-K,K]\big] &\ge 6^n \mathbb{P}[A_{n,\delta_0}(v)]\Big(
\prod_{k\in \widetilde{I}_{4K,n}} \mathbb{P}\big[M_{n-k-1}\le 0\big]\mathbb{P}\big[M_{n-k-1}\ge 0\big]
\Big[\farc{1}{2}\int_{ (-\delta_0,\delta_0)^c} q(x)\dx\Big]^{2}\Big)\\
&\times \Big[\prod_{k\in \widetilde{I}_{4K,n}^c} \mathbb{P}[M_{n-k-1}\le -3K]\mathbb{P}[M_{n-k-1}\ge 3K]\Big]\\
&\ge 6^n\mathbb{P}[A_{n,\delta_0}(v)] \Big(\frac{1}{16}\Big[\int_{ (-\delta_0,\delta_0)^c} q(x)\dx\Big]^2\Big)^{|\widetilde{I}_{4K,n}|}
\Big(\frac{1}{2}-\frac{1}{200}\Big)^{2|\widetilde{I}_{4K,n}^c|}.
\enbal
\ee
In the last step, we used the symmetry of $q$ for $k\in\widetilde{I}_{4K,n}$, while for $k\in \widetilde{I}_{4K,n}^c$
we used (\ref{sep308}) 
together with the definition (\ref{sep310}) of $\widetilde{I}_{4K,n}^c$. 

Now, assume that 
\be\label{sep320}
\liminf_{n\to \infty} \frac{|\widetilde{I}_{4K,n}|}{n} = 0.
\ee
\begin{lem}\label{lem-sep502}
Assume that (\ref{sep320}) holds. Then, there exists $\eta_0>0$ so that for all $\eta\in (0,\eta_0)$ 
there is~$C_\eta>0$ such that for $K\ge C_\eta$ and all $n\in\N$ we have 
\be\label{sep321}
\liminf_{n\to \infty}\farc{\mathbb{P}[A_{n,\delta_0}(v)]}{ (1-\eta)^n} \ge 1.
\ee
\end{lem}
We postpone the proof of this lemma for the moment.
Fix $\eta>0$ sufficiently small and fix $\varepsilon>0$ such that  
\be\label{sep322}
\farc{3}{2.5} (1-\eta)\Big (1-\farc{1}{100}\Big)^2\cdot \Big(\int_{ (-\delta_0,\delta_0)^c} q(x)\dx\Big)^{2\varepsilon} >1.
\ee
Note that (\ref{sep322}), together (\ref{sep320}) and $|\widetilde{I}_{4K,n}^c|\le n$,  implies that
\be
\bal
1&\ge \liminf_{n\to \infty}\mathbb{P}[M_n\in[-K,K]]\ge \liminf_{n\to \infty} 6^n(1-\eta)^n\farc{1}{(16)^{2\eps n}}
\Big(\int_{ (-\delta_0,\delta_0)^c} q(x)\dx\Big)^{\varepsilon n} \Big(\frac{1}{2}-\frac{1}{200}\Big)^{2n}
\\
&\ge \liminf_{n\to \infty}\Big(\frac{3}{2.5}\Big)^n (1-\eta)^n 
\Big(\int_{ (-\delta_0,\delta_0)^c} q(x)\dx\Big)^{\varepsilon n} \Big(1-\frac{1}{100}\Big)^{2n}>1,
\enbal
\ee
yielding a contradiction. Thus (\ref{sep320}) can not hold. This gives 
\be\label{sep323}
\liminf\limits_{n\to\infty} \frac{|I_{4K}\cap \{0,\dots, n-1\}|}{n} = \liminf\limits_{n\to\infty} \frac{|\widetilde{I}_{4K,n}|}{n}>0.
\ee
Since, as we have observed at the beginning of the proof, $(M_n)_{n\in I_{4K}}$ is tight, (\ref{sep323}) yields the claim of Theorem~\ref{Theo:SubseqTight}.~$\Box$

\subsubsection*{The proof of Lemma~\ref{lem-sep502}}

Let us define
\be\label{sep716}
\pi_{n,\delta_0, m,C} := \inf_{I\subseteq\{1,\dots, n\}, |I|\le m} \mathbb{P}\left[|S_k|\le C\ \text{for all}\ k\le n\ \text{and}\ |S_k|\le \delta_0\ \text{for all}\ k\in I \right].
\ee
We will show that for all $c<1$ there are $C_0>0$, $\varepsilon_1>0$ such that for all $n\in\N$ big enough we have
\begin{equation}
\pi_{n,\delta_0, \lfloor \varepsilon_1 n \rfloor, C_0}\ge c^n. \label{eq:pigr}
\end{equation}
First, we show how \eqref{sep320} and \eqref{eq:pigr} imply \eqref{sep321}. We use the definition (\ref{sep518}) of $A_{n,\delta_0}$ to write,
%
%
for $\varepsilon_1>0$ arbitrary, $n$ big enough, depending on $\varepsilon_1$, and $v$ such that $|v|=n$:
\be\label{sep520}
\bal
\mathbb{P}[A_{n,\delta_0}(v)] &=\Pm\Big[|S_{v_k}-S_v|\le K \hbox{ for all $k\le n$} \hbox{ and }   |S_v-S_{v_k}|\le \delta_0
\hbox{ for all } k\in \widetilde{I}_{4K,n}\Big]
\\
&= \mathbb{P}\left[|S_{n-k}|\le K\ \text{for all}\ k\le n\ \text{and}\ |S_{n-k}|\le \delta_0\ \text{for all}\ k\in \widetilde{I}_{4K,n}\right]\\
&=\mathbb{P}\left[|S_k|\le K\ \text{for all}\ k\le n\ \text{and}\ |S_k|\le \delta_0\ \text{for all}\ k\in n-\widetilde{I}_{4K,n}\right]
\ge \pi_{n,\delta_0,\lfloor \varepsilon_1 n\rfloor, K}.
\enbal
\ee
The second equality above used the symmetry of $q$ and the equivalence
\[
(S_{v_k}-S_v)_{k\le n} \stackrel{d}{=} (-S_{n-k})_{k\le n},
\]
while the last inequality in (\ref{sep520}) used assumption \eqref{sep320}. 

Now, \eqref{eq:pigr} and (\ref{sep520}) imply that there are $\varepsilon_1$ and $C_0$, which depend on $\eta$, such that for $K\ge C_0$ and~$n$ big enough we have
\[
\mathbb{P}[A_{n,\delta_0}(v)] \ge (1-\eta)^n,
\]
which implies \eqref{sep321}. 

The idea of the proof of \eqref{eq:pigr} is to use the Markov property at all times in $I$ and also bound from below 
the probability that between these times the random walk remains in $[-C,C]$ and ends in $[-\delta_0,\delta_0]$.
Thus, we define for $N\in\N$
\be\label{sep522} 
p_{0,\delta_0,C,N}:=\inf_{x\in [-\delta_0,\delta_0]}\mathbb{P}\left[\forall_{k\le N} |x+S_k|\le C, |x+S_N|\le \delta_0 \right].
\ee
We will choose $L_C>0$, split the time interval $[0,N]$ into intervals of length $L_C$, and 
force~$|x+S_k|\le C/2$ at the end of these pieces. We will also use the Markov property at the start of each of these intervals.
We will need a slightly different calculation for the last piece and will also need to deal with the case~$N<L_C$. 
It will be helpful to use the following notation
%
\be 
\bal
p_{1,\delta_0} &:= \inf_{C\ge2\delta_0}\inf_{x\in [-C/2,C/2]} \mathbb{P}\left[\forall_{k\le L_C} |x+S_k|\le C, |x+S_{L_C}|\le C/2 \right],\\
p_{2,\delta_0,C} &:= \inf_{k\in \{1,\dots, L_C\}} \inf_{x\in [-\delta_0,\delta_0]} \mathbb{P}\left[ \forall_{j\le k} |x+S_j|\le C, |x+S_k|\le \delta_0\right],\\
p_{3,\delta_0, C} &:= \inf_{k\in \{L_C+1,\dots, 2L_C\}} \inf_{x\in [-C/2,C/2]}\mathbb{P}\left[\forall_{j\le k} |x+S_j|\le C, |x+S_k|\le \delta_0\right].
\enbal
\ee
Note that if $N\le L_C$ then 
\begin{equation} \label{sep524}
p_{0,\delta_0, C,N} \ge  p_{2,\delta_0,C},
\end{equation}
while if $L_C<N\le 2L_C$ then 
\begin{equation} \label{sep525}
p_{0,\delta_0, C,N} \ge  p_{3,\delta_0, C},
\end{equation}
and if $N>2L_C$ then
\begin{equation} \label{sep526}
p_{0,\delta_0,C,N} \ge p_{1,\delta_0}^{{N}/{L_C}}p_{3,\delta_0,C}.
\end{equation}
Together, (\ref{sep524}), (\ref{sep525}) and \eqref{sep526} imply
\begin{equation} \label{eq:p0Nzwischen}
p_{0,\delta_0, C,N} \ge p_{1,\delta_0}^{ {N}/{L_C}} p_{2,\delta_0, C}p_{3,\delta_0, C}.
\end{equation}
To make use of \eqref{eq:p0Nzwischen} we need to prove that all three factors are strictly positive. For $p_{2,\delta_0,C}$ we have 
\be\label{sep512}
\bal
p_{2,\delta_0,C} &\ge \inf_{k\in \{1,\dots,  L_C\}} \inf_{x\in [-\delta_0,\delta_0]} \mathbb{P}\left[ \forall_{j\le k} |x+S_j|\le \delta_0, |x+S_k|\le \delta_0\right] 
\ge \Big(\inf_{x\in[-\delta_0,\delta_0]} \int_{-\delta_0}^{\delta_0} q(y-x)\dy\Big)^{L_C}>0,
\enbal
\ee
due to the choice of $\delta_0$ in (\ref{sep506}). 

Next, we prove that there is a $C_1\ge 0$ such that for $C\ge C_1$ we can choose $L_C$ such that 
\be\label{sep702}
\hbox{$p_{3,\delta_0,C}>0$ and $L_C\to \infty$ as $C\to\infty$. }
\ee
First, by symmetry it is enough to prove that 
\be\label{sep712}
\inf_{k\in \{L_C+1,\dots, 2L_C\}} \inf_{x\in [-C/2,0]}\mathbb{P}\left[\forall_{j\le k} |x+S_j|\le C, |x+S_k|\le \delta_0\right] >0.
\ee
The idea is to force the random walk to drift towards 0 with increments smaller than $\delta_0$, such that it can't skip over the interval $[-\delta_0,\delta_0]$.
Once the random walk hits the interval, we can use the second condition in \eqref{sep506} to force the random walk to stay inside $[-\delta_0,\delta_0]$ until 
the time $2L_C$, at a cost smaller than $\gamma^{2L_C}$, with some $\gamma>0$. However, as we do not require $q$ to have mass near $0$ we cannot 
force it to have a small increment in every step. Instead, we use the symmetry of $q$ to  force the two-step increment~$S_{k+2}-S_k$ to be small.
To this end, we claim that, as $q$ is continuous and symmetric, there are intervals $I_1\subseteq [0,\infty)$,~$I_2\subseteq (-\infty,0]$ and $k_q\in \N$ such that
\be \label{eq:p0gr0}
p_0 := \min\Big[\int_{I_1} q(y)\dy,~\int_{I_2} q(y)\dy\Big]>0,
\ee
and 
\be
\hbox{$z_1+z_2 \in [\delta_0/(2k_q),\delta_0/k_q]$, for all $z_1\in I_1$, $z_2\in I_2$.}
\ee
Moreover, we can chose $I_1$ and $I_2$ to be of the form
\be\label{sep706}
I_1 :=\Big[\delta_0 \Big(\frac{l_q}{k_q}+\frac{3}{4k_q}\Big), \delta_0 \frac{l_q+1}{k_q}\Big],
~~~ I_2 := \Big[-\delta_0 \Big(\frac{l_q}{k_q}-\frac{1}{4k_q}\Big), -\delta_0 \frac{l_q}{k_q}\Big],
\ee 
with some $l_q\in\Nm$. 
We set 
\be
C_1 := 2\delta_0\frac{l_q+1}{k_q}.
\ee
Then, for all $C\ge C_1$, we have
\be\label{sep708}
\hbox{$x+z \in [-C/2,C/2]$, for all $x\in[-C/2,0]$ and $z\in I_1$.}
\ee
We set
\be\label{sep710}
L_C := \left\lfloor \frac{C}{2}\frac{2 k_q}{\delta_0}\right\rfloor.
\ee

Next, consider the stopping time 
\[
\tau_1^x := \inf\{k\in\N : x+S_k \in [-\delta_0,\delta_0]\},
\]
and, for $T\in\Nm$, the event
\[
B_{x,T} := \{\forall_{j\le T : k\in 2\N+1} S_j \in I_1,~~  \forall_{j\le T : k\in 2\N} S_j \in I_2\}.
\]
We note that, by the choice of $I_1$ and $I_2$ in (\ref{sep706}), for $x\in [-C/2,0]$ on $B_{x,\tau_1^x}$ we have  
\be
\hbox{$x+S_k \in [-C/2,C/2]$ for all $k\le \tau_1^x$.}  
\ee
Moreover, the choice (\ref{sep710}) of $L_C$ implies that $\tau_1^x\le L_C$.
It follows that
\be\label{sep514}
\bal
&\inf_{k\in \{L_C+1,\dots, 2L_C\}} \inf_{x\in [-C/2,0]}\mathbb{P}\left[\forall_{j\le k} |x+S_j|\le C, |x+S_k|\le \delta_0\right]\\
& \ge \inf_{k\in \{L_C+1,\dots, 2L_C\}} \inf_{x\in [-C/2,0]} \mathbb{P}\left[ B_{x,\tau_1^x},\ \forall_{j \in \{\tau_1^x+1,\dots,k\}} |x+S_j| \le \delta_0\right]\\
&\ge \inf_{k\in \{L_C+1,\dots, L_C\}} \inf_{x\in [-C/2,0]} \inf_{T\in \{0,\dots, L_C\}} \left( \mathbb{P}[B_{x,T}] 
\cdot \inf_{z\in [-\delta_0,\delta_0]} \mathbb{P}\left[ \forall_{j\le k-T} |z+S_{j}| \le \delta_0 \right] \right)\\
&\ge  p_0^{L_C}\cdot \left(\inf_{z\in [-\delta_0,\delta_0]} \int_{-\delta_0}^{\delta_0} q(y-z)\dy\right)^{2L_C}>0,
\enbal
\ee
with $p_0$ as in (\ref{eq:p0gr0}). 
The last step above used the second condition on $\delta_0$ in \eqref{sep506}. This proves (\ref{sep712}). Thus, (\ref{sep702}) is also proved.

 Finally, the inequality $p_{1,\delta_0}>0$ can be seen using a path-wise version of the CLT and the definition of $L_C$ in \eqref{sep710}. Overall, we have proven that the right side of (\ref{eq:p0Nzwischen}) is positive. 

Now, fix $J_n$ with $|J_n| \le \lfloor \varepsilon_1 n\rfloor+2$. If needed, we can add $0$ and $n$ to $J_n$. 
Let $\{x_k\}$ be an ordered enumeration of $J_n$. Using the Markov property at the times $x_k$ we see that 
\be \label{sep718}
\bal
&\mathbb{P}\left[|S_k|\le C\ \text{for all}\ k\le n\ \text{and}\ |S_k|\le \delta_0\ \text{for all}\ k\in J_n \right] \ge \prod_{k=0}^{|J_n|-1} p_{0,\delta_0,C,x_{k+1}-x_k}\\
&
~~~~~~~~~~~~~~\qquad
{\ge} \prod_{k=0}^{|J_n|-1} \Big(p_{1,\delta_0}^{({x_{k+1}-x_k})/{L_C}} p_{2,\delta_0, C}p_{3,\delta_0, C} \Big)
\ge p_{1,\delta_0}^{{n}/{L_C}}p_{2,\delta_0,C}^{\varepsilon_1 n+2}p_{3,\delta_0,C}^{\varepsilon_1 n+2}.
\enbal
\ee
We used (\ref{eq:p0Nzwischen}) in the second inequality above. 
Going back to the definition (\ref{sep716}) of $\pi_{n,\delta_0, m,C}$, we deduce from (\ref{sep718}) that 
\be\label{sep722}
\pi_{n,\delta_0, \lfloor \varepsilon_1 n\rfloor, C}\ge p_{1,\delta_0}^{{n}/L_C} p_{2,\delta_0,C}^{\varepsilon_1 n+2}p_{3,\delta_0,C}^{\varepsilon_1 n+2}.
\ee

Finally, observe that, given any  $c<1$, we can take $C$ sufficiently large, so that 
\be\label{sep720}
p_{1,\delta_0}^{{1}/{L_C}}\ge  c^{1/3}.
\ee
Next, fix $\varepsilon_1>0$ such that $(p_{2,\delta_0,C_0}p_{3,\delta_0,C_0})^{\varepsilon_1}\ge c^{1/3}$ and $n$ sufficiently large, so that
\be\label{sep721}
p_{2,\delta_0, C_0}^2p_{3,\delta_0,C_0}^2 \ge c^{(1/3)n}.
\ee
Together, (\ref{sep722})--(\ref{sep721}) imply 
\be\label{sep723}
\pi_{n,\delta_0,\lfloor \varepsilon_1 n\rfloor, C_0} \ge c^n.
\ee
This proves (\ref{eq:pigr}) and finishes the proof of Lemma~\ref{lem-sep502}.~$\Box$ 

\section{Tightness in the single zero bistable case with analytic means} \label{Sec:TightnessAnMeans}

In this section, we consider, by analytic means, random threshold 
voting models for which the nonlinearity~$f(u)$ defined by~(\ref{eq:gThresh})--(\ref{aug2802})
has exactly one zero $\vtheta\in(0,1)$. In addition, we assume that
\be\label{sep726}
f'(0)<0,
~~f'(1)<0.
\ee  
In particular, it follows that
\be\label{sep727}
\hbox{$f(x)<0$ for $x\in(0,\vtheta)$, $f(x)>0$ for $x\in(\vtheta,1)$.}
\ee
This is the bistable case: the zeroes $x=0$ and $x=1$  of $f(x)$ are stable and $x=\vtheta$ is an unstable zero. We will extend the nonlinearity $f(x)$
and the recursion polynomial $g(x)$ outside of $[0,1]$ by setting
\be\label{sep728}
\bal
&f(x)=f'(0)x,~~g(x)=x+f'(0)x,~~\hbox{ for $x<0$},\\
&f(x)=f'(1)(x-1),~~g(x)=x+f'(1)(x-1),~\hbox{for $x>1$.}
\enbal
\ee
Let us comment that since $g(x)$ corresponds to a random threshold voting model,
Proposition~\ref{Theo:NonLinThresh} implies that $g(x)=x+f(x)$ is increasing on $(0,1)$. 
Hence, in addition to (\ref{sep727}), we must have 
\be
-1\le f'(0),f'(1)<0.
\ee
It follows that the extension of $g(x)$ is non-decreasing on all of $\Rm$. A standard example of
such nonlinearity is the binary-ternary voting model described in Section~\ref{subsec-bt}, with $p<1/2$. 

In addition to the standing assumptions (\ref{aug3025}) on $q$, we  
assume that $q\in C^1(\R)$.
Under these conditions we will first prove the following theorem.
\begin{thm}\label{Theo:Tightness}
Let $(p_d,\zeta_{k,d})_{d\le d_0, k\in\{1,\dots, d\}}$ be a random threshold model such that the nonlinearity $f(u)$ satisfies (\ref{sep726})-(\ref{sep727}). 
Then, the sequence $(M_n-\mathrm{med}(M_n))_{n\in\N}$ is tight.
\end{thm}

The next step is we show in Theorem~\ref{thm:sep1102-thm}, 
using the result of Theorem~\ref{Theo:Tightness}, 
 that $\mathrm{med}(M_n)$ itself has the asymptotics
\be\label{nov2702}
\hbox{med}(M_n)=n\ell+x_0+o(1),~~\hbox{as $n\to+\infty$.}
\ee
Here, $\ell$ is the speed of a unique traveling wave constructed in Proposition~\ref{Lem:ExTravWav}
below. 
We also show in this theorem that the distribution $\Pm(M_n>x)$ converges to a shift of the traveling wave,
strengthening the tightness claim of Theorem~\ref{Theo:Tightness}.
Note that the  conclusion in (\ref{nov2702}) 
differs from the classical maximum of branching random walks setup, were $(n\ell-\mathrm{med}(F_{M_n}))$ is of 
the order $\log n$.

The proof of Theorem \ref{Theo:Tightness} is divided into two steps. First, we use~\cite{yagisita2010existence} to show that there exists 
a traveling wave solution to the recursion \eqref{eq:RecEqThresh}.  In the second step, we use a discrete in time version of the Fife-McLeod technique~\cite{FifeMcL}
to prove that $F_{M_n}$ can be bound between a super-solution and a sub-solution to \eqref{eq:RecEqThresh}, which are constructed by perturbing 
the traveling wave solution. This, in particular, shows the uniqueness of the traveling wave speed.
Here, the bistable assumptions~(\ref{sep726})--(\ref{sep727}) on $f(x)$ are essential.  

\subsection{Existence of a traveling wave}

A traveling wave is a solution to  \eqref{eq:RecEqThresh}
\begin{equation}\label{sep729}
w_{{n+1}} = g(q\ast w_{n}), 
\end{equation}
of the form
\be\label{sep728bis}
w_n(x)=\vphi(x-n\ell),
\ee
with some $\ell\in\Rm$. We say that $\ell$ is the speed of the wave and $\vphi$ is its profile. Equivalently,  the traveling wave
is a solution to 
\begin{equation}\label{eq:trwvcond}
\varphi  =g(q_\ell\ast \varphi),  
\end{equation}
with $q_l(x) := q(x+\ell)$, together with the boundary conditions
\be\label{sep808}
\vphi(-\infty)=0,~~\vphi(+\infty)=1.
\ee
Indeed, if $\vphi(x)$ satisfies (\ref{eq:trwvcond}) then 
$w_n(x) := \varphi(x-n\ell)$ satisfies
\begin{equation}
\begin{aligned}
w_{n+1}(x) &= \varphi(x-(n+1)\ell) = g(q_\ell\ast \varphi(\cdot - (n+1)\ell))(x) = g\Big( \int_{\R} q_\ell(y) \varphi(x-y-(n+1)\ell)\dy\Big)\\
&= g\Big(\int_{\R} q(y+\ell)\varphi(x-y-(n+1)\ell)\dy\Big) = g\left(\int_{\R} q(y) \varphi(x-y-n\ell)\dy\right) = g(q\ast w_n)(x),
\end{aligned}
\end{equation}
which is \eqref{sep729}.

We will use the following comparison principle for \eqref{sep729}.
As a notation, 
we let $\cal M$ be the set of monotone non-decreasing and left continuous functions $w(x)$ on $\Rm$ such that the
limits 
\be
w_\pm=\lim_{x\to\pm\infty}w(x)
\ee
exist and are finite. For an interval $I=[w_-,w_+]$, we will denote by ${\cal M}_I$ the set of functions in ${\cal M}$ with the
corresponding left and right limits. 
\begin{prop}\label{Lem:Comparison}
Suppose that the sequence $\{w_n\}_{n\in\N}\subseteq \mathcal{M} $ is a solution to~\eqref{eq:RecEqThresh}.\\
(i) If $\{\overline{w}_n\}_{n\in\N} \subseteq \mathcal{M}$ satisfies
\begin{equation}\label{eq:SuperSol}
\overline{w}_{n+1}(x) \ge g(q\ast \overline{w}_n)(x),~~\hbox{  for all $n\ge 0$ and $x\in\R$,}
\end{equation}
and $\overline{w}_0(x)\ge w_0(x)$ for all $x\in \R$, then  $\overline{w}_n(x)\ge w_n(x)$ for all $n\ge 0$ and~$x\in\R$. 
\\
(ii) If $\{\underline{w}\}_{n\in\N}\subseteq \mathcal{M} $ satisfies
\begin{equation} \label{eq:SubSol}
\underline{w}_{n+1}(x) \le g(q\ast \underline{w}_n)(x),~~\hbox{ for all $n\ge 0$ and $x\in\R$},
\end{equation}
and $\underline{w}_0(x)\le w_0(x)$ for all $x\in\R$, then $\underline{w}_n(x)\le w_n(x)$ for all $n\ge0$ and $x\in\R$.
\end{prop}
The proof of this proposition is immediate, once one 
recalls that by Proposition~\ref{Theo:NonLinThresh}, the function $g(x)$ is increasing on $(0,1)$ (and its extension in \eqref{sep728}
continues to be increasing outside that interval). 
This is the main reason that we can only handle monotone recursion polynomials~$g$ in this section, that is,
recursion equations corresponding to random threshold models.

The next proposition gives the existence of a traveling wave. 
\begin{prop}[Existence of a traveling wave] \label{Lem:ExTravWav}
There exist   $\ell\in\R$ and  a non-decreasing $\varphi\in C^1(\R)$ that satisfy~(\ref{eq:trwvcond})--(\ref{sep808}). 
\end{prop}
We will see later that both the speed $\ell$ and the travelling wave $\vphi$ are unique (up to a shift of the latter).
\begin{proof}[Proof of Proposition \ref{Lem:ExTravWav}]
%
The claim of Proposition~\ref{Lem:ExTravWav} follows from 
Corollary 5 of~\cite{yagisita2010existence}. Let us briefly explain the details. 
Setting 
\be
Q_0[u](x)  := g(q\ast u)(x),
\ee
we can write the traveling wave equation (\ref{eq:trwvcond}) as 
\be\label{sep730}
\vphi(x)=Q_0[\varphi](x+\ell) .
\ee  
The aforementioned corollary establishes the existence of a non-decreasing solution $\vphi(x)$
to (\ref{sep730}) that satisfies the boundary conditions (\ref{sep808}) under the following assumptions (Hypotheses 2 and 3
in~\cite{yagisita2010existence}):\\ 
(i) The map $Q_0$ is continuous  with respect to locally uniform convergence. That is, if 
$\{u_k\}_{k\in\N}\subseteq \mathcal{M}_{[0,1]}$ converges to $u\in\mathcal{M}_{[0,1]}$ 
uniformly on every bounded interval, the sequence $\{Q_0[u_k]\}_{k\in\N}$ converges to $Q_0[u]$ almost everywhere.\\
(ii) The map $Q_0$ is order-preserving.  \\
(iii) The map $Q_0$ is
translation invariant.\\
(iv) The map $Q_0$ is bistable, in the sense that  that there is $\alpha\in (0,1)$ with $Q_0[\alpha] = \alpha$, $Q_0[\gamma]<\gamma$ 
for all~$0<\gamma<\alpha$ and $\gamma<Q_0[\gamma]$ for all~$\alpha<\gamma<1$.\\
(v) If there are two constants $\ell_-, \ell_+\in \R$ and non-decreasing functions $\varphi_-$ and $\vphi_+$ such that 
\be\label{sep732}
\hbox{$(Q_0[\varphi_-])(x+\ell_-) = \varphi_-(x)$, $\varphi_-(-\infty) = 0$, $\varphi_-(+\infty) = \vartheta$,}
\ee
and
\be\label{sep733}
\hbox{ $(Q_0[\varphi_+])(x+\ell_+) = \varphi_+(x)$, $\varphi_+(-\infty) = \vartheta$ and $\varphi_+(+\infty) =1$,}
\ee
then 
\be\label{sep731}
\ell_->\ell_+.
\ee
This means that any traveling wave solution to \eqref{sep729} connecting 
$0$ to $\vartheta$ travels faster to the right than a traveling wave solution connecting $\vartheta$ to $1$.


It is straightforward to verify that assumptions (i)--(iv) above are satisfied here. 
In particular, continuity and translation invariance in assumptions (i) and (iii) follow immediately from the definition of~$Q_0$ and our assumptions on $q(x)$.
The order preserving property (ii) is a consequence of the comparison principle in Proposition~\ref{Lem:Comparison}. 
The bistable assumption (\ref{sep727}) on the nonlinearity $f(x)$ implies assumption~(iv) above.

The last step is to verify assumption (v) above on the speed comparison. Let $\vphi_-$ and $\vphi_+$ be, respectively, solutions to (\ref{sep732}) and (\ref{sep733}). 
We first consider $\vphi_-$ and write (\ref{sep732}) as 
\begin{equation}\label{sep734}
\varphi_-(x-\ell_-) =g(u(x)),~~u(x)=\int q(y)\vphi_-(x-y)\dy
\end{equation}
that can be written as 
\be
\vphi_-(x-\ell_-)-\vphi_-(x)=\int q(y)[\vphi_-(x-y)-\vphi_-(x)]\dy+g(u(x))-u(x),
\ee
or, equivalently, as 
\be
\vphi_-(x-\ell_-)-\vphi_-(x)=\int q(y)[\vphi_-(x-y)-\vphi_-(x)]\dy+f(u(x)).
\ee
Integrating in $x\in[-M,M]$ with some $M\gg 1$ gives
\be\label{sep804}
\bal
\int_{-M}^M[\vphi_-(x-\ell_-)-\vphi_-(x)]\dx&=\int_{-M}^M\int_\Rm q(y)[\vphi_-(x-y)-\vphi_-(x)]\dy\dx+\int_{-M}^Mf(u(x))\dx\\
&=
\int_\Rm q(y)G_M(y)\dy+\int_{-M}^Mf(u(x))\dx,
\enbal
\ee
with
\be\label{sep802}
G_M(y)=\int_{-M}^M [\vphi_-(x-y)-\vphi_-(x)]\dx=\int_{-M-y}^{-M}\vphi_-(x)\dx-\int_{M-y}^M\vphi_-(x)\dx.
\ee
Passing to the limit $M\to+\infty$ in (\ref{sep802}) using the boundary conditions in (\ref{sep732}) gives
\be
\lim_{M\to+\infty} G_M(y)=-y\vtheta.
\ee
Similarly, passing to the limit in the left side of (\ref{sep804}) gives
\be
\lim_{M\to+\infty}\int_{-M}^M[\vphi_-(x-\ell_-)-\vphi_-(x)]\dx=-\ell_-\vtheta.
\ee
In addition, the boundary conditions in (\ref{sep732}) 
and (\ref{sep727}) imply that 
\be
\lim_{M\to+\infty}\int_{-M}^M f(u(x))\dx<0.
\ee
Thus, passing to the limit $M\to+\infty$ in (\ref{sep804}) gives
\be
-\ell_-\vtheta< -E_q\vtheta,
\ee
with
\be
E_q=\int_\Rm yq(y)\dy.
\ee
We conclude that
\be
\ell_-> E_q.
\ee
A completely analogous argument shows that
\be
\ell_+<E_q.
\ee
Now, (\ref{sep731}) follows. Therefore, assumption (v) also holds, and Corollary 5 of~\cite{yagisita2010existence} can be applied. 
This finishes the proof of Proposition \ref{Lem:ExTravWav}. 
\end{proof}

\subsection{Basic properties of a traveling wave}

We now prove some basic properties of any traveling wave that will be needed in the proof of Theorem~\ref{Theo:Tightness} as well as Theorem~\ref{thm:sep1102-thm} below. 
First, we get a bound on the traveling wave speed $\ell$. 
\begin{lem}\label{Lem:WaveSpeed}
If $\varphi(x)$ and $\ell\in\Rm$ satisfy \eqref{eq:trwvcond}--(\ref{sep808}), 
then $\ell\in \left(\min\supp(q),\max\supp(q)\right)$.
\end{lem}
\begin{proof}
Let $M_{n,\varphi}$ be the  outcome of the threshold voting model associated to the recursion polynomial~$g$ from \eqref{eq:trwvcond}, 
where the starting location of the underlying branching random walk is distributed according to $\varphi$. 
Similarly to \eqref{eq:RecEqThresh},   $M_{n,\varphi}$ solves
\be \label{sep810}
\bal
&F_{M_{n+1,\varphi}}(x) = g\left(q\ast F_{M_n,\varphi}\right)(x),\\
&F_{M_{0,\varphi}}(x) = \varphi(x).
\enbal
\ee
As $\vphi(x)$ is a traveling wave, we know that the solution to (\ref{sep810}) is 
\be
F_{M_{n,\varphi}}(x) = \varphi(x-n\ell).
\ee
Fix $c\in \R$ with $\varphi(c) = 1/2$, so that 
\[
\mathbb{P}[M_{n,\varphi}\le c+n\ell]  = \varphi(c+n\ell-n\ell) = \varphi(c) = \frac{1}{2}.
\]
Since 
the distribution of $M_{n,\varphi}$ is continuous we also have
\be
\mathbb{P}[M_{n,\varphi}\ge c+n\ell]=\frac{1}{2}.
\ee

Assume, for the sake of contradiction, that 
\be\label{sep814}
\ell\ge \max \supp(q).
\ee
Let us choose $r\in(1/2,1)$ such that there is a unique $q_r$ so that $\vphi(q_r)=r$. 
Also, let $X_0\sim \varphi$ be independent of the underlying BRW. We have  
\be\label{sep812}
\bal
\frac{1}{2} &= \mathbb{P}[M_{n,\varphi}\ge c+n\ell] = \mathbb{P}[M_n+X_0\ge c+n\ell]
\le \mathbb{P}[X_0 > q_{r}]+\mathbb{P}[X_0\le q_{r}, M_n\ge c+n\ell-X_0]\\
&\le 1-r 
+\mathbb{P}\left[M_n\ge c+n\ell-q_{r}\right] \le 1-r 
+\mathbb{P}\left[\exists_{v : |v| = n} S_v \ge c+n\ell-q_{r} \right] \\
&\le 1-r 
+d_0^n\mathbb{P}[S_n\ge c+n\ell-q_{r}],
\enbal
\ee
where 
\[
S_n = \sum_{k=1}^n X_k,~~(X_k)_{k\in\N}\hbox{ i.i.d.\@ and $X_1\sim q$},
\]
and $d_0$ is the maximal number of children one particle can have, so that the total number of particles  in generation $n$ 
is bounded by $d_0^n$.

Since $q$ is non-atomic and $\ell\ge \max\supp(q)$ there is  
$\ell_{d_0}<\ell$ such that for $n$ big enough we have
\be\label{sep825}
\mathbb{P}[S_n\ge n\ell_{d_0}] \le (d_0+1)^{-n}.
\ee
To see that \eqref{sep825} holds, fix $1/2>\eta>0$ and $\delta_{\eta}>0$ 
such that \[
\mathbb{P}[X_1\ge \max \supp(q)-2\delta_{\eta}]\le \eta,
\]
and set 
\[
I_{1,\eta} := (-\infty, \max\supp(q)-2 \delta_{\eta}],~~
I_{2,\eta} := [\max\supp(q)-2\delta_{\eta}, \infty).
\]
Such $\delta_{\eta}$ exists since $q$ has no atoms, 
so that~$X_1$ has a continuous distribution function.
Then, we have
\be
\bal
S_n &\le |\{k\le n : X_k \in I_{1,\eta}\}|\cdot (\max\supp(q)-2\delta_{\eta})+|\{k\le n : X_k \in I_{2,\eta}\}|\cdot \max\supp(q) \\
&= n\max\supp(q)-2\delta_{\eta}|I_{1,\eta}|.
\enbal
\ee
This implies 
\be \label{eq:SnZuBin}
\mathbb{P}[S_n \ge n(\max\supp(q)-\delta_{\eta})] \le \mathbb{P}\left[|I_{1,\eta}| \le \frac{n}{2}\right] \le \mathrm{Bin}\left(n,1-\eta\right)[0,n/2] = \mathrm{Bin}\left(n,\eta\right)[n/2,n].
\ee
By Corollary 2.2.19 and Exercise 2.2.23~(b) in \cite{LDPDZ},~\eqref{eq:SnZuBin} implies  
\be \label{eq:LDP}
\limsup_{n\to\infty} \frac{1}{n}\log\left(\mathbb{P}[S_n\ge n(\supp(q)-\delta_{\eta})\right) \le -\inf_{x\ge 1/2} \left(x\log\left(\frac{x}{\eta}\right)+(1-x)\log\left(\frac{1-x}{1-\eta}\right) \right)
\ee
Since 
\[
\lim\limits_{\eta\to 0} \inf_{x\ge 1/2} 
\left(x\log\left(\frac{x}{\eta}\right)+(1-x)\log\left(\frac{1-x}{1-\eta}\right)\right) = \infty,
\]
the inequality \eqref{eq:LDP} implies \eqref{sep825}.

We now choose $\ell_{d_0}$ as in (\ref{sep825}). Then, 
for $n$ large enough, we have 
\[
c+n\ell-q_{r}>n\ell_{d_0}.
\]
Thus, (\ref{sep812}) yields that, for $n$ big enough, we have 
\[
\frac{1}{2}\le 1-r 
+d_0^n\mathbb{P}[S_n\ge n \ell_{d_0}] \le 1-r +\Big(\frac{d_0}{d_0+1}\Big)^n.
\]
Passing to the limit $n\to+\infty$ gives a contradiction since $r>1/2$. Thus, (\ref{sep814}) can not hold, whence
\[
\ell<\max\supp(q).
\] 
An analogous argument starting with  
\[
\frac{1}{2} = \mathbb{P}[M_{n,\varphi}\le c+n\ell]
\]
yields that $\ell>\min\supp(q)$, which finishes the proof of Lemma \ref{Lem:WaveSpeed}.
\end{proof}

The next lemma shows that the traveling wave profile has no critical points.
\begin{lem}\label{Lem:StrictIncr}
Any traveling wave solution $\varphi$ to \eqref{eq:trwvcond} has $\varphi(x)\in (0,1)$ and $\varphi'(x)>0$ for all $x\in \R$.
\end{lem}
\begin{proof}
We will prove that $\varphi(x)<1$ for all $x\in\R$, the proof that $\varphi(x)>0$ for all $x\in\R$ is analogous. 
Assume, for the sake of contradiction, that there is  some $x\in\R$ with $\varphi(x) = 1$ and consider
\be\label{sep818}
x_0:=\min\{x:~\vphi(x)=1\}.
\ee 
Since $\varphi$ is a solution to \eqref{eq:trwvcond}, we have 
\[
1 = \varphi(x_0) = g\left(\int_{\R} q(y+\ell)\varphi(x_0-y)\dy\right).
\]
Since $g(x) = 1$ iff $x = 1$, we deduce 
\be
1 = \int_{\R} q(y)\varphi(x_0-y+\ell)\dy,
\ee 
which, in turn, implies that for all $y\in \supp(q)$ we have 
\begin{equation}
\varphi(x_0-y+\ell) = 1.
\end{equation} 
However, by Lemma \ref{Lem:WaveSpeed} we know that $\ell<\max\supp(q)$. Thus, there is some $y\in\supp(q)$ with 
\[
x_0-y+\ell<x_0.
\]
This is a contradiction to the definition (\ref{sep818}) of $x_0$.

Next, we prove that $\varphi'(x)>0$ for all $x\in\R$. Differentiating \eqref{eq:trwvcond} yields 
\begin{align}
\varphi'(x) = g'((q_\ell\ast \varphi)(x))\cdot \int_{\R} q(y+\ell)\varphi'(x-y)\dy = g'((q_\ell\ast \varphi)(x)) \cdot \int_{\R} q(y)\varphi'(x-y+\ell)\dy. \label{eq:TrWvDeriv}
\end{align}
Assume, for the sake of a contradiction, that there is some $x_0\in \R$ with 
\be\label{sep820}
\varphi'(x_0)=0.
\ee
Since $\varphi(x)\in (0,1)$ for all $x\in\R$ 
and $g'(u)>0$ for all $u\in (0,1)$, (\ref{eq:TrWvDeriv}) implies  
\be\label{sep821}
\int_{\R} q(y)\varphi'(x_0+\ell-y)\dy = 0.
\ee
Let $I_q$ be an interval on which $q(x)$ is strictly positive. 
Since $\varphi$ is non-decreasing, it follows from (\ref{sep821}) that 
\be
\hbox{$\varphi'(x) = 0$ for all $x\in x_0+\ell-I_q$.}
\ee
Iterating this argument, we conclude that 
\be
\hbox{$\varphi'(x) = 0$ for all $x\in x_0+n\ell-n\cdot I_q$ for all $n\ge 1$.}
\ee
Thus, there is an arbitrarily long interval $I\subseteq \R$ on which $\varphi(x)\equiv z$ is constant. As we have shown that~$\vphi(x)$ takes values in $(0,1)$,
we have $z\in(0,1)$. In addition, if $I$ is sufficiently long, $z$ must be a solution to 
\be
g(z)=z.
\ee
It follows that $z=\vtheta$. As such intervals are arbitrarily long and $\vphi(x)$ is non-decreasing, this leads to a contradiction to the boundary conditions for $\vphi(x)$. 
\end{proof}

\subsection{The proof of Theorem~\ref{Theo:Tightness}}

The proof of Theorem~\ref{Theo:Tightness} relies on the following trapping of the solution $F_{M_n}$ to (\ref{sep729})
between two perturbations of the traveling wave solution.
\begin{lem}\label{Lem:MainIngr}
There exists an increasing bounded sequence $\xi_n^+$ and an decreasing bounded sequence $\xi_n^-$ and constants $\beta_0^+, \beta_0^->0$, $\delta_0^+, \delta_0^->0$ such that
\begin{equation}
\overline{w}_n(x) = \varphi(x-n\ell+\xi_n^+)+\beta_0^+ e^{-\delta_0^+n}
\end{equation}
satisfies \eqref{eq:SuperSol} and
\begin{equation}
\underline{w}_n(x) = \varphi(x-n\ell+\xi_n^-)-\beta_0^- e^{-\delta_0^-n}
\end{equation}
satisfies \eqref{eq:SubSol}. Furthermore, we can choose $\xi_0^+$ arbitrarily large and $\xi_0^-$ arbitrarily small without changing~$\beta_0$, $\delta_0$. 
\end{lem}
Here, by convention, we extend $g(u)$ outside of $[0,1]$ as in \eqref{sep728}. 
Note that the extension is still an increasing function, so that the comparison principle
in Proposition \ref{Lem:Comparison} still applies. 
Before we prove Lemma \ref{Lem:MainIngr}, we show how it implies Theorem \ref{Theo:Tightness}.

\begin{proof}[Proof of Theorem \ref{Theo:Tightness} assuming Lemma \ref{Lem:MainIngr}]
We first show that by choosing $\xi_0^+$, $\xi_0^-$ in Lemma \ref{Lem:MainIngr} appropriately we can assure that for all $n\in\N_0$ and $x\in\R$ we have 
\begin{equation}\label{eq:Sandwich}
\underline{w}_n(x)\le F_{M_n}(x)\le \overline{w}_n(x).  
\end{equation}
Using Proposition \ref{Lem:Comparison} and \eqref{eq:RecEqThresh} reduces the proof of 
\eqref{eq:Sandwich} to showing that we can choose $\xi_0^+$, $\xi_0^-$ such that 
\be
\underline{w}_0(x) \le\one(x\ge 0) \le \overline{w}_0(x),
\ee 
which is easy to arrange because $\varphi(-\infty) = 0$ and $\varphi(+\infty) = 1$.


The definitions of $\overline{w}_n$ and $\underline{w}_n$ and
\eqref{eq:Sandwich} imply
tightness of $(M_n-n\ell)_{n\in\N}$, which in particular implies that 
\[
\sup_n (n\ell-\mathrm{med}(F_{M_n})) <\infty,
\]
and thus that $(M_n-\mathrm{med}(M_n))_{n\in\N}$ is tight as well.
\end{proof}

\subsubsection*{The proof of Lemma \ref{Lem:MainIngr}}

We write
\[
\beta^+_n := \beta^+_0e^{-\delta_0^+ n},
\]
with $\beta^+_0$, $\delta_0^+$ to be chosen later on. A function of the form 
\[
\overline{w}_n(x) = \varphi(x-n\ell+\xi^+_n)+\beta^+_n
\]
satisfies \eqref{eq:SuperSol} if
\begin{align*}
0&\le N_n(x) := \overline{w}_{n+1}-g(q\ast \overline{w}_n) = \varphi(x-(n+1)\ell+\xi^+_{n+1})+\beta^+_{n+1}
-g\left(q\ast (\varphi(x -n\ell+\xi^+_n)+\beta^+_n)\right)\\
&\stackrel{\eqref{eq:trwvcond}}{=} g\big(q_\ell\ast \varphi(x -(n+1)\ell+\xi^+_{n+1})\big)+\beta^+_{n+1}-
g\Big(\int_{\R} q(y)\varphi(x-n\ell+\xi^+_n-y)\dy+\beta^+_n\Big)\\
&= g\Big(\int_{\R} q(y+\ell)\varphi(x-y-(n+1)\ell+\xi^+_{n+1})\dy\Big)+\beta^+_{n+1}
-g\Big(\int_{\R} q(y)\varphi(x-n\ell+\xi^+_n-y)\dy+\beta^+_n\Big)\\
&= g\Big(\int_{\R} q(y)\varphi(x-y-n\ell+\xi^+_{n+1})\dy\Big)+\beta^+_{n+1}
-g\Big(\int_{\R} q(y)\varphi(x-n\ell+\xi^+_n-y)\dy+\beta^+_n\Big).
\end{align*}
We set
\[
\zeta_n^+ := x-n\ell+\xi^+_n,
\]
and consider the regions $|\zeta_n^+|\ge R_0$, $|\zeta_n^+|\le R_0$ separately.
Here, $R_0>0$ will be chosen sufficiently large later on.

\paragraph{The exterior region $|\zeta_n^+|\ge R_0$.} Let us set
\be\label{sep1104}
I_n := \int_{\R} q(y)\varphi(\zeta_n^+-y)\dy.
\ee
Since the sequence $\xi^+_n$ will be chosen non-decreasing,
and $\varphi$ is non-decreasing as well, we have 
\begin{align}
N_n &\ge g\left(\int_{\R} q(y)\varphi(x-y-n\ell+\xi^+_n)\dy\right)
+\beta^+_{n+1}-g\left(\int_{\R} q(y)\varphi(x-y-n\ell+\xi^+_n)\dy+\beta^+_n\right)\notag\\
&= g\left(\int_{\R} q(y)\varphi(\zeta_n^+-y)\dy\right)+\beta^+_{n+1}-g\left(\int_{\R} q(y)\varphi(\zeta_n^+-y)\dy+\beta^+_n\right)\notag\\
& = (g(I_n)-I_n)+\beta^+_{n+1}+I_n - \left(g(I_n+\beta^+_n)-(I_n+\beta^+_n)\right)-(I_n+\beta^+_n) \notag\\
&=\beta^+_{n+1}-\beta^+_n+ f(I_n)-f(I_n+\beta^+_n). \label{eq:znbig}
\end{align}
Let us recall that by (\ref{sep726}) we have 
$f'(0)<0$ and $f'(1)<0$. Furthermore, by increasing $R_0$, we can make both $I_n$ and $1-I_n$
be arbitrarily close to zero in the region $|\zeta^+_n|\ge R_0$.
In particular, we can choose $R_0$ large, $\gamma_0>0$ small, such that 
for $|\xi^+_n|\ge R_0$ and $0\le \beta^+_n\le \gamma_0$ we have
\be\label{sep1102}
f(I_n)-f(I_n+\beta^+_n)\ge \eta_0\beta^+_n.
\ee
Here, $\eta_0>0$ can be chosen, for example, as 
\[
\eta_0 = \frac{1}{2}\min\left(|f'(0)|, |f'(1)|\right)>0.
\]
Using (\ref{sep1102})  in \eqref{eq:znbig} shows that $N_n\ge 0$ 
in the region $|\zeta_n^+|\ge R_0$ if 
\[
\beta^+_{n+1}-\beta^+_n+\eta_0\beta^+_n\ge0.
\]
This is true if we take $\beta^+_0\in(0,\gamma_0)$ and set 
\be
\beta^+_n = \beta^+_0(1-\eta_0)^n = \beta^+_0 e^{-\log\left( (1-\eta_0)^{-1}\right)n},
\ee 
with $\delta_0^+ =  -\log(1-\eta_0)$.

\paragraph{The interior region $|\zeta^+_n|\le R_0$.}
Let $C_g$ be the Lipschitz constant of $g(u)$ and write
\begin{align}
N_n &=g\left(\int_{\R} q(y)\varphi(\zeta_n^++\xi^+_{n+1}-\xi^+_n-y)\dy\right)+\beta^+_{n+1}-g\left(\int_{\R} q(y)\varphi(\zeta_n^+-y)\dy+\beta^+_n\right)\notag\\
&\ge g\left(\int_{\R} q(y)\varphi(\zeta_n^++\xi^+_{n+1}-\xi^+_n-y)\dy\right)+\beta^+_{n+1}-g\left(\int_{\R} q(y)\varphi(\zeta_n^+-y)\dy\right)-C_g\beta^+_n \notag\\
&= g(I_n+E_n)-g(I_n)+\beta^+_{n+1}-C_g\beta^+_n, \label{eq:Nnznsm}
\end{align}
with $I_n$ as in (\ref{sep1104}) and 
\[
E_n := \int_{\R} q(y)(\varphi(\zeta_n^++\xi^+_{n+1}-\xi^+_n-y)-\varphi(\zeta_n^+-y)) \dy.
\]
We assume now that, in addition to $|\zeta_n^+|\le R_0$, we have 
\be\label{sep1106}
|\xi^+_{n+1}-\xi^+_n|\le 100.
\ee 
Note that, since $\varphi$ is monotone and $\xi_n^+$ is increasing in $n$,
we have
\be\label{eq:Enlb}
E_n \ge \int_{-R_0}^{R_0} q(y)(\varphi(\zeta_n^++\xi^+_{n+1}-\xi^+_n-y)-\varphi(\zeta_n^+-y)) 
\dy\ge \eta_1 (\xi^+_{n+1}-\xi^+_n).
\ee
Here, we have set 
\[
\eta_1 := \int_{-R_0}^{R_0}q(y)\dy\cdot \inf_{x\in [-2R_0-100,2R_0+100]} \varphi'(x).
\] 
Observe that, after potentially increasing $R_0$, so that $q([-R_0,R_0])>0$, 
and using Lemma \ref{Lem:StrictIncr}, we know that 
\be
\eta_1>0.
\ee

Next, we note that, since 
\be
I_n= (q\ast\varphi)(\zeta_n^+),~~
I_n+E_n = (q\ast \varphi)(\zeta_n^++\xi^+_{n+1}-\xi^+_n), 
\ee
and $|\zeta_n^+|\le R_0$, 
there is a $\delta_1>0$ such that
\begin{equation}
\delta_1<I_n<I_n+E_n<1-\delta_1. \label{eq:In01}
\end{equation}
We set
\be
\eta_2 := \min_{u\in (\delta_1,1-\delta_1)} g'(u)>0.
\ee

Combining \eqref{eq:Enlb}, \eqref{eq:In01} and \eqref{eq:Nnznsm} gives
\begin{equation}
N_n\ge \eta_2E_n+\beta^+_{n+1}-C_g\beta^+_n \ge \eta_1\eta_2(\xi^+_{n+1}-\xi^+_n)-C_g\beta^+_n. \label{eq:Nn2ndlast}
\end{equation}
Thus, to ensure that $N_n\ge 0$ in the region $|\zeta_n^+|\le R_0$, 
it is enough to take
\be\label{sep1118}
\xi^+_{n+1} = \xi^+_n+\frac{C_g}{\eta_1\eta_2}\beta^+_n 
= \xi^+_n+\frac{C_g}{\eta_1\eta_2}\beta^+_0(1-\eta_0)^n,
\ee
so that
\be\label{sep1121}
\xi_n^+=\xi_0^++\frac{C_g}{\eta_1\eta_2}\beta^+_0\sum_{k=0}^{n-1}(1-\eta_0)^k
=\xi_0^+ +K_g\beta_0^+(1-\gamma_0^n\big),
\ee
with appropriately defined $K_g>0$ and $\gamma_0=1-\eta_0$. 
This sequence is increasing and bounded. Moreover, for $\beta^+_0$ small enough, but, 
importantly, independent of~$\xi^+_0$, we have (\ref{sep1106})  as well. 

Thus, we have shown that we can find $\beta^+_0, \delta_0^+>0$ such that for all $\xi_0^+\ge 0$ there is an increasing bounded sequence $(\xi_n^+)_{n\in\N_0}$ such that 
\be
\overline{w}_n(x) = \varphi(x-nl+\xi_n^+)+\beta^+_0e^{-\delta_0^+n}
\ee 
satisfies \eqref{eq:SuperSol}. The corresponding construction of $\underline{w}_n(x)$
that satisfies (\ref{eq:SubSol}) is very similar. We only mention that the sequence $\xi_n^-$ can be chosen as
\be\label{sep1123}
\xi_n^-=\xi_0^- -K_g\beta_0^-(1-\gamma_0^n) .
\ee
This finishes the proof of Lemma~\ref{Lem:MainIngr}.~$\Box$

Let us finish this section with the following corollary of Lemma~\ref{Lem:MainIngr} and its proof.
\begin{cor}\label{lem-sep1104}
There exist $K>0$, $\delta_0>0$,  and $r_0>0$ with the following property. Suppose that~$w_n(x)$ is a solution
to the recursion 
\begin{equation}\label{nov829}
w_{{n+1}} = g(q\ast w_{n}), 
\end{equation}
with an initial condition $w_0(x)$ that satisfies
\be\label{nov830}
\vphi(x+\xi_0^-)-\beta_0\le w_0(x)\le \vphi(x+\xi_0^+)+\beta_0,~~\hbox{for all   $x\in\Rm$}.
\ee
Then, if $0\le \beta_0\le \delta_0$, we have
\be\label{sep1124}
\vphi(x-n\ell+\xi_n^-)-\beta_n\le w_n(x)\le \vphi(x-n\ell+\xi_n^+)+\beta_n,
~~\hbox{for all $n\ge 1$ and $x\in\Rm$},
\ee
with
\be\label{nov810}
\beta_n= \beta_0\exp(-r_0n),
\ee
and
\be\label{sep1134}
\xi_-^n=\xi_0^- -K\beta_0\big(1-\exp(-r_0n)\big),~~
\xi_n^+=\xi_0^+ +K\beta_0\big(1-\exp(-r_0n)\big).
\ee 
\end{cor}
We remark that Corollary \ref{lem-sep1104} implies that the speed $\ell$ 
in Proposition \ref{Lem:ExTravWav} is unique 
but not yet the uniqueness of the traveling wave profile $\vphi(x)$.

\subsection{The long time convergence to a traveling wave}

Our goal here is to prove the following theorem. 
\begin{thm}\label{thm:sep1102-thm}
Under the assumptions of Theorem \ref{Theo:Tightness}, 
let $\vphi$ be a travelling wave as in Proposition \ref{Lem:ExTravWav}, shifted so that $\vphi(0)=1/2$. Let $\ell$ be its speed.
Then, there exists $x_0\in\Rm$ such that
\be\label{sep1112}
\med (M_n)=n\ell+x_0+o(1),~~\hbox{as $n\to+\infty$,}
\ee
and
\be\label{sep1304}
|\Pm(M_n<x)-\vphi(x-n\ell-x_0)|\to 0.
\ee
\end{thm}
In the proof of this theorem, it will be more convenient to work not with $w_n(x)=\Pm(M_n<x)$ but its translation in the moving
frame of the traveling wave
\be
u_n(x)=w_n(x+n\ell).
\ee
The function $u_n(x)$ satisfies the recursion (\ref{eq:trwvcond})
\be\label{nov802}
\bal
&u_{n+1}(x)=g((q_\ell\ast u_n)(x)),
&u_0(x)=\one(x\ge 0),
\enbal
\ee
for which $\vphi(x)$ is a fixed point:
\be
\vphi(x)=g((q_\ell\ast \vphi)(x)).
\ee
Here, as in (\ref{eq:trwvcond}), we have set 
\be
q_\ell(x)=q(x+\ell).
\ee
More generally, we will assume that the initial condition
$u_0(x)$ satisfies 
\be\label{nov1310}
\vphi(x+\xi_0^-)-\beta_0\le u_0(x)\le \vphi(x+\xi_0^+)+\beta_0,~~\hbox{for all   $x\in\Rm$},
\ee
with some $\xi_0^\pm\in\Rm$ and  $0\le \beta_0\le \delta_0$.
Corollary~\ref{lem-sep1104} implies that then $u_n(x)$
obeys the uniform bounds
\be\label{nov804}
\vphi(x+\xi_n^-)-\beta_n\le u_n(x)\le \vphi(x+\xi_n^+)+\beta_n,
~~\hbox{for all $n\ge 1$ and $x\in\Rm$},
\ee
with $\beta_n\to 0$ according to (\ref{nov810}), and uniformly bounded $\xi_n^\pm$:
\be\label{nov812}
|\xi_n^\pm|\le K,~~\hbox{ for all $n\ge 0$.}
\ee

Together with the a priori regularity estimates on $u_n(x)$, this implies,
in particular, that the iterates~$\{u_n(\cdot)\}_n$ lie in a compact subset of $C(\Rm)$. We will denote
by ${\cal Z}[u_0]$ the $\omega$-limit set of $\{u_n(\cdot)\}_n$. It consists of
all functions~$\zeta_n(x)$, defined for $n\in\Zm$ and $x\in\Rm$, such that there is a sequence~$n_k\in\Nm$ (independent of $n$) so that  
$n_k\to+\infty$ as $k\to+\infty$, and 
\be\label{nov806}
u_{n+n_k}(x)\to \zeta_n(x),~~\hbox{as $k\to+\infty$.}
\ee
The limit in (\ref{nov806}) is uniform on $x\in\Rm$ and finite sets of $n\in\Zm$.
Note that any such limit~$\zeta_n(x)$ is a global in time solution to (\ref{nov802}), defined for all
$x\in\Rm$ and $n\in\Zm$:
\be\label{nov831}
\bal
&\zeta_{n+1}(x)=g((q_\ell\ast \zeta_n)(x)),
\enbal
\ee
with the initial condition 
\be\label{nov1302}
\zeta_0(x)=\lim_{k\to+\infty}u_{n_k}(x).
\ee
An important point is that the solution to (\ref{nov831}) with the initial condition as in (\ref{nov1302}) is defined both for~$n\ge 0$ and~$n\le 0$.
Let us stress that the set ${\cal Z}[u_0]$ depends on the choice of the initial condition $u_0$ for~(\ref{nov802}). Another helpful observation is that if $\zeta\in\cZ[u_0]$, 
then $\cZ[\zeta_k]\subseteq\cZ[u_0]$, for any $k\in\Zm$ fixed. 

An immediate consequence of the bounds in (\ref{nov804}), as well as (\ref{nov810}) and (\ref{nov812}) is that there
exist $\bar\xi_\pm$ so that any element~$\zeta\in{\cZ}[u_0]$ satisfies
\be\label{nov814}
\vphi(x+\bar\xi_-)\le \zeta_n(x)\le \vphi(x+\bar\xi_+),~~\hbox{ for all $n\in\Zm$ and $x\in\Rm$.}
\ee
Our goal is to show that $\cZ[u_0]$ contains exactly one element and that element is a traveling wave.
The key step is the following.
\begin{prop}\label{prop-nov802}
The $\omega$-limit set $\cZ[u_0]$ contains a traveling wave $\vphi(x+\bar\xi)$, with some $\bar\xi\in\Rm$. 
\end{prop}
Proposition~\ref{prop-nov802}, together with the stability estimates
(\ref{sep1124})-(\ref{sep1134}) in Corollary~\ref{lem-sep1104},  implies
immediately that there is exactly one traveling wave in the $\omega$-limit set of $\{u_n\}_n$, and that this traveling wave is the only element
of $\cZ[u_0]$, finishing the proof of Theorem~\ref{thm:sep1102-thm}. 

We first prove the following lemma.
\begin{lem}\label{lem-nov1302}
Let $\zeta\in\cZ[u_0]$ and suppose that for some
$\xi\in\Rm$ we have
\be\label{nov1320}
\zeta_n(x)\le\vphi(x+\xi),~~\hbox{ for all $n\in\Zm$ and $x\in\Rm$.}
\ee
Assume, in addition, that there exist $n_0\in\Zm$ and
$y_0\in\Rm$ so that 
\be\label{nov1306}
\zeta_{n_0}(y_0)=\vphi(y_0+\xi).
\ee
Then, we have
\be\label{nov1321}
\zeta_n(x)=\vphi(x+\xi),~~\hbox{ for all $n\in\Zm$ and $x\in\Rm$.}
\ee
\end{lem}
{\bf Proof.} We may assume without loss of generality that $y_0=0$. 
Our assumptions on $q$ and the result of Lemma~\ref{Lem:WaveSpeed} 
that $\ell$ is srtrictly inside the support of $q$ implies
that there exist two intervals $I_-=[y_-,y_+]$ and $I_+=[x_-,x_+]$ with $y_{\pm}<0$,
 $x_{\pm}>0$ and
$q_\ell>0$ on $I_-\cup I_+$.
Since both $\zeta_n(x)$ and $\vphi(x+\xi)$ are solutions to the recursion (\ref{nov802}),
one obtains 
that if both (\ref{nov1320})
and (\ref{nov1306}) hold, then
\begin{equation}
  \label{eq-nov15}
  \zeta_{n_0-n}(x)=\varphi(x+\xi), \quad \mbox{\rm for $x\in \cup_{k\leq n}
  (kI_-+(n-k)I_+)$}.
\end{equation}
Further,
since
there exist $k,n$ so that $kI_-+(n-k)I_+$ contains an interval around $0$, one 
deduces
(by taking multiples of such $n$)
that for each $R>0$ there exists $n_R$ so that
\be\label{nov1322}
\zeta_{n}(x)=\vphi(x+\xi),~\hbox{for all $|x|\le R$ and $n\le n_0-n_R$}.
\ee
Recall that there exist $\delta_1>0$ and $\delta_2>0$ so that 
\be\label{nov1323}
f'(u)<-\delta_1,~~\hbox{ for $u\in[0,\delta_2]$ and $u\in[1-\delta_2,1]$.}
\ee
To use this stability in the tails, 
we will choose $R_0>0$, 
so that for all $\zeta\in \cZ[u_0]$,
\be\label{nov1324}
\bal
&0\le\zeta_n(x),~\vphi(x+\xi)\le\delta_2,~~\hbox{ for all $x\le -R_0$ and $n\in\Zm$},\\
&1-\delta_2 \le\zeta_n(x),~\vphi(x+\xi)\le 1,~~\hbox{ for all $x\ge R_0$ and $n\in\Zm$}.
\enbal
\ee
This is possible due to the estimates in (\ref{nov814}). As $q(x)$ is compactly supported and has
mass equal to one, it follows
that there exists~$M>0$ so that for all $\zeta\in \cZ[u_0]$,
\be\label{nov1330}
\bal
&0\le q\ast\zeta_n(x),~q\ast\vphi(x+\xi)\le\delta_2,~~\hbox{ for all $x\le -R_0-M$ and $n\in\Zm$},\\
&1-\delta_2 \le q\ast\zeta_n(x),~q\ast\vphi(x+\xi)\le 1,~~\hbox{ for all $x\ge R_0+M$ and $n\in\Zm$}.
\enbal
\ee

Let us consider the difference
\be\label{nov1325}
y_n(x):=\vphi(x+\xi)-\zeta_n(x).
\ee
Note that, because of (\ref{nov1322}), given any $M>0$, we know that there exists $m_R$ so that
\be\label{nov1326}
y_n(x)=0,~~\hbox{ for all $n\le n_0-m_{R}$ and $|x|\le R_0+M$.}
\ee
In the region $|x|\ge R_0+M$, the function $y_n(x)$ satisfies an equation of the form 
\be\label{nov1327}
y_{n+1}(x)=(q_\ell\ast y_n)(x)+a_n(x)(q_\ell\ast y_n)(x),
\ee
with, recalling $\delta_1$ from \eqref{nov1323},
\be\label{nov1328}
a_n(x)=\left\{\begin{array}{ll}
  \farc{f(q\ast z)(x)-f(q\ast\zeta_n)(x)}{(q\ast z)(x)-(q\ast\zeta_n)(x)},& q\ast z(x)\neq q\ast \zeta_n(x)\\
  -\delta_1, & \mbox{\rm else},
\end{array}\right.
\ee
and
\be
z(x):=\vphi(x + \xi).
\ee
We know from  (\ref{nov1330}) that for $|x|\ge R_0+M$ the arguments of $f(\cdot)$ in
(\ref{nov1328}) are sufficiently close to $0$ on the left and $1$ on the right, and
we have 
\be\label{nov1329}
a_n(x)\le -\delta_1,~~\hbox{for $|x|\ge R_0+M$ and $n\in\Zm$.}
\ee
Therefore, if we choose $M>0$ sufficiently large, depending only on the support
of $q(\cdot)$, we will have
\be\label{nov1331a}
y_{n+1}(x)\le (q_\ell\ast y_n)(x)-\delta_1(q_\ell\ast y_n)(x),~~\hbox{ for all $|x|\ge R_0+M$ and $n\in\Zm$}.
\ee
 
Next, observe that (\ref{nov1326}) implies that for any $n\le n_0-n_R$ the non-negative
function $y_n(x)$ attains its maximum 
\be
Y_n=\max_{x\in\Rm}y_n(x),
\ee
at some point $x_n$ such that $|x_n|\ge R_0+M$:
\be
Y_n=y_n(x_n).
\ee
It follows now from (\ref{nov1331a}) that 
\be\label{nov1335}
Y_{n}\le (1-\delta_1)^mY_{n-m},~~\hbox{for all $n\le n_0-n_R$.}
\ee
Letting $m\to+\infty$ with $n\le n_0-n_R$ fixed, we deduce that $Y_n=0$, which, in turn, implies
that $y_n(x)\equiv 0$ for all $n\le n_0-n_R$. This, of course, implies that $y_n(x)\equiv 0$ for all $n\in\Zm$,
which is (\ref{nov1321}).~$\Box$


{\bf Proof of Proposition~\ref{prop-nov802}.}  
Consider an arbitrary element $\zeta\in\cZ[u_0]$ and set 
\be\label{nov816}
\bar\xi_{sm}[\zeta]=\inf\{\xi:~\hbox{ $\zeta_n(x)\le\vphi(x+\xi)$,
for all $x\in\Rm$ and $n\in\Zm$}\}.
\ee
It follows from (\ref{nov814}) that 
\be\label{nov1308}
\bar\xi_-\le\bar\xi_{sm}[\zeta]\le \bar\xi_+,~~\hbox{ for all $\zeta\in\cZ[u_0]$.}
\ee
Note that, in particular, we have
\be\label{nov818}
\zeta_n(x)\le\vphi(x+\bar\xi_{sm}[\zeta]),~\hbox{for all $x\in\Rm$ and $n\in\Zm$}.
\ee
Let us consider
\be\label{nov823}
\bar\xi_{sm}=\inf_{\zeta\in\cZ[u_0]}\bar\xi_{sm}[\zeta].
\ee
As the set $\cZ[u_0]$ is compact, if we take a sequence $\zeta^{(k)}\in\cZ[u_0]$ such 
that $\bar\xi_{sm}[\zeta^{(k)}]\to \bar\xi_{sm}$,
then, after passing to the limit~$\zeta^{(k)}\to\bar\zeta$, possibly along a subsequence, we will 
find $\bar\zeta\in\cZ[u_0]$ such that
\be\label{nov824}
\bar\xi_{sm}[\bar\zeta]=\bar\xi_{sm}.
\ee
As in (\ref{nov818}), we will still have
\be\label{nov1314}
\bar\zeta_n(x)\le\vphi(x+\bar\xi_{sm}[\zeta]),~\hbox{for all $x\in\Rm$ and $n\in\Zm$}.
\ee
We deduce from Lemma~\ref{lem-nov1302} that either we have 
\be\label{nov821}
\bar\zeta_n(x)=\vphi(x+\bar\xi_{sm}[\zeta]),~\hbox{for all $x\in\Rm$ and $n\in\Zm$},
\ee
and we are done, or  
\be\label{nov1316}
\bar\zeta_{n}(x)<\vphi(x+\bar\xi_{sm}[\zeta]),~\hbox{for all $x\in \Rm$ and $n\in\Zm$}.
\ee
Suppose that (\ref{nov1316}) holds and let $\cZ[\bar\zeta_0]\subset\cZ[u_0]$ be  
the $\omega$-limit set of $\bar\zeta_0(x)$.  
We claim that either $\cZ[\bar\zeta_0]$ contains
a shift of a traveling wave, in which case so does $\cZ[u_0]$, and we are done, or not only
(\ref{nov1316}) holds but also for any $R>0$ 
there exist $\delta_R>0$ and $n_R\in\Zm$ so that 
\be\label{nov822}
\vphi(x+\bar\xi_{sm})-\bar\zeta_n(x)\ge \delta_R>0,~~\hbox{ for all $n\ge n_R$ and $|x|\le R$.}
\ee
Indeed, otherwise there would exist a sequence $n_k\to+\infty$ and $x_k\in[-R_0,R_0]$ such that 
\be\label{nov1331}
\vphi(x_k+\bar\xi_{sm})-\bar\zeta_{n_k}(x_k)\to0,~~\hbox{ as $k\to+\infty$.}
\ee
Therefore, possibly after further extracting a sub-sequence, we would find an element $\eta\in\cZ[\bar\zeta_0]$
and a point~$y\in[-R_0,R_0]$ such that
\be\label{nov1332}
\eta_{n}(x)\le\vphi(x+\bar\xi_{sm}[\zeta]),~\hbox{for all $x\in \Rm$ and $n\in\Zm$},
\ee
and
\be
\eta_0(y)=\vphi(y+\bar\xi_{sm}[\zeta]).
\ee
Lemma~\ref{lem-nov1302} would then imply that 
\be
\eta_n(x)\equiv\vphi(x+\bar\xi_{sm}[\zeta]).
\ee
Therefore, the set $\cZ[u_0]$ would contain a traveling wave.

To finish the proof, we will show that if (\ref{nov822}) holds, 
then there is an element~$\tilde\zeta\in\cZ[u_0]$ such that 
\be\label{nov827}
\bar\xi_{sm}[\tilde\zeta]<\bar\xi_{sm},
\ee
which will be a contradiction to the definition of $\bar\xi_{sm}$. 
Let us suppose that $\bar\zeta_n(x)$ is a solution to (\ref{nov802})
such that (\ref{nov822}) holds for all $R>0$ and, in addition, we know that  
\be\label{nov827a}
\vphi(x+\bar\xi_{sm})-\bar\zeta_n(x)\ge  0,~~\hbox{ for all $n\in\Zm$ and $x\in\Rm$.}
\ee
Once again, after extracting a subsequence and passing to the limit, we will find an element 
$\eta\in\cZ[\bar\zeta_0]$ such that the restriction $n\ge n_R$ in (\ref{nov822}) can be removed:
\be\label{nov1334}
\vphi(x+\bar\xi_{sm})-\eta_n(x)\ge \delta_R>0,~~\hbox{ for all $n\in \Zm$ and $|x|\le R$.}
\ee
%
We argue as at the end of the proof of Lemma~\ref{lem-nov1302}. 
Note that, because of (\ref{nov1334}), given any $M>0$, we can take $\gamma>0$ sufficiently small, so that, with $R_0$ as in \eqref{nov1324}, 
\be\label{nov833}
y_n(x):=\vphi(x+\bar\xi_{sm}-\gamma)-\eta_n(x)\ge \farc{\delta_{R_0}}{2}>0,~~\hbox{ for all $n\in\Zm$ and $|x|\le R_0+M$.}
\ee
As in the aforementioned proof, in the region $|x|\ge R_0+M$, 
the function $y_n(x)$ satisfies an equation of the form 
\be\label{nov834}
y_{n+1}(x)=(q_\ell\ast y_n)(x)+a_n(x)(q_\ell\ast y_n)(x),
\ee
with $a_n(x)$ as in \eqref{nov1328} but with $\zeta_n$ replaced by $\eta_n$ and 
\be
z(x):=\vphi(x + \bar\xi_{sm}-\gamma).
\ee
If $R_0$ and $M$ are chosen as in \eqref{nov1324} and \eqref{nov1330},
we have that
\be\label{nov837}
a_n(x)\le -\delta_1,~~\hbox{for $|x|\ge R_0$ and $n\in\Zm$.}
\ee
In addition, outside of this region, we have, using (\ref{nov833}) 
\be\label{nov838}
y_{n}(x)>\farc{\delta_{R_0}}{2},~~\hbox{ for all $n\in\Zm$ and $|x|\le R_0+M$. }
\ee
Moreover, at any "initial" time $m$ we have
\be\label{nov839}
y_m(x)>-C\gamma,~~\hbox{ for $|x|\ge R_0+M$.}
\ee
In particular, using \eqref{nov833} and setting $y_m^*=\min_{x\in \Rm} y_m(x)\wedge 0$,
we have that $y_{m}^*\geq (1-\delta_1) y_m^*$. Hence,
\be\label{nov840}
y_{m+k}(x)\ge -C\gamma e^{-\delta_1k},~~\hbox{ for $|x|\ge R_0+M$.}
\ee
As the starting time $m$ is arbitrary, it follows that actually $y_n(x)\ge 0$ for all $x\in\Rm$ and $n\in\Zm$. 
Therefore, we have
\be\label{nov841}
\vphi(x+\bar\xi_{sm}-\gamma)-\bar\zeta_n(x)\ge  0,~~\hbox{ for all $n\in\Zm$ and $x\in\Rm$.}
\ee
As $\gamma>0$, this contradicts the definition of $\bar\xi_{sm}$, 
finishing the proof of Proposition~\ref{prop-nov802}. The proof of Theorem~\ref{thm:sep1102-thm} 
is complete as well.~$\Box$

\bibliographystyle{plain}

\end{document}